\documentclass[preprint,12pt]{elsarticle}

\usepackage{amssymb}
\usepackage{amsthm}
\usepackage{amsmath}
\usepackage{amsfonts}
\usepackage{mathrsfs}
\usepackage{upgreek}
\usepackage{mathtools}
\usepackage{bm}

\usepackage{xspace}
\usepackage{tabularx}
\usepackage[caption=false]{subfig}

\usepackage{algorithm}
\usepackage{algorithmicx,setspace}
\usepackage{algpseudocode}
\usepackage[dvipsnames]{xcolor}
\usepackage{url}

\usepackage{tikz-cd}
\usetikzlibrary{cd}
\usepackage{adjustbox}
\usepackage{graphicx}
\usepackage{glossaries}

\usepackage[margin=2.5cm]{geometry}

\usepackage{setspace} 

\journal{}

\DeclareMathOperator\supp{supp}

\newacronym{fem}{FEM}{finite element method}
\newacronym{fd}{FD}{finite differences}
\newacronym{fft}{FFT}{fast Fourier transform}
\newacronym[longplural={degrees of freedom}]{dof}{DOF}{degree of freedom}

\newcommand{\D}[1]{\mathbb{#1}}

\newcommand{\HTran}{{\mathsf{H}}}
\newcommand{\Transpose}{{\mkern-1.5mu\mathsf{T}}}

\newcommand{\N}{N}
\newcommand{\Nsubtext}[1]{{N_\mathrm{#1}}}

\newcommand{\NI}{\Nsubtext{I}}
\newcommand{\Np}{\Nsubtext{p}}
\newcommand{\Nn}{\Nsubtext{n}}

\newcommand{\NQ}{\Nsubtext{Q}}

\usepackage{boites, boites_exemples}
\newcounter{till_comment_counter}

\newcommand{\puc}{\mathcal{Y}}
\newcommand{\intd}[1]{\,{\mathrm d}#1}
\newcommand{\MB}[1]{\boldsymbol{\mathsf{#1}}} 

\newcommand{\vek}[1]{\mathchoice{\displaystyle\boldsymbol{#1}}
{\textstyle\boldsymbol{#1}}{\scriptstyle\boldsymbol{#1}}
{\scriptscriptstyle\boldsymbol{#1}}}

\newcommand{\mat}[1]{\mathchoice{\displaystyle\mathbf{#1}}
{\textstyle\mathbf{#1}}{\scriptstyle\mathbf{#1}}
{\scriptscriptstyle\mathbf{#1}}}

\newcommand{\norm}[1]{\left\lVert#1\right\rVert}

\newcommand{\alp}{\ensuremath{\alpha}}

\newcommand{\flux}{\ensuremath{\vek{q}}}
\newcommand{\stress}{\ensuremath{\vek{\sigma}}}

\newcommand{\strain}{\ensuremath{\vek{\varepsilon}}}
\newcommand{\macrostrain}{\ensuremath{\vek{e}}}

\newcommand{\grad}{\ensuremath{\nabla}}
\newcommand{\symgrad}{\ensuremath{\grad_{\text{s}}}}
\newcommand{\symgradM}{\ensuremath{\vek{\partial}}}

\newcommand{\dispcomponent}{\ensuremath{u}}
\newcommand{\disp}{\vek{\dispcomponent}}

\newcommand{\perdispcomponent}{\Tilde{\dispcomponent}}
\newcommand{\perdisp}{\vek{\Tilde{\dispcomponent}}}
\newcommand{\Dperdisp}{\MB{\Tilde{\dispcomponent}}}

\newcommand{\testf}{\Tilde{\ensuremath{v}}}
\newcommand{\perdisptest}{\vek{\testf}}

\newcommand{\qweights}{\ensuremath{W}}

\newcommand{\materialsymbol}{C}
\newcommand{\material}{\mat{\materialsymbol}}
\newcommand{\Dmaterial}{\MB{\materialsymbol}}
\newcommand{\rmaterial}{\mat{\materialsymbol}^\text{ref}}

\newcommand{\Drmaterial}{\MB{\materialsymbol}^\text{ref}}

\newcommand{\testspace}{\ensuremath{\mathcal{V}}}

\newcommand{\identity}{I}

\newcommand{\FourierTrans}{\mathcal{F}}
\newcommand{\iFourierTrans}{\mathcal{F}^{-1}}

\newcommand{\Doper}{\mathcal{D}}
\newcommand{\DoperTrans}{\mathcal{D}^{\Transpose}}
\newcommand{\Dmat}{\MB{D}}

\newcommand{\precMatrix}{\MB{M}}

\newcommand{\Kref}{\MB{K}_{(i)}^{\text{ref}}}

\newcommand{\KrefF}{\widehat{\MB{K}}{}_{(i)}^{\textnormal{ref}}}

\newcommand{\Krefab}{                \MB{K}   _{(i)}^{\text{ref}}{}_{\bar{\alp}\bar{\beta}}}
\newcommand{\hatKrefab}{    \widehat{\MB{K}}{}_{(i)}^{\text{ref}}{}_{\bar{\alp}\bar{\beta}}}
\newcommand{\KrefFab}{      \widehat{\MB{K}}{}_{(i)}^{\text{ref}}{}_{\bar{\alp}\bar{\beta}}}

\newcommand{\Krefcomp}[2] {         {\MB{K} _{(i)}^{\text{ref}}}{}_{#1#2}}
\newcommand{\KrefFcomp}[2]{ \widehat{\MB{K}}{}_{(i)}^{\text{ref}}{}_{#1#2}}

\newcommand{\Krefinv} {{(\MB{K}_{(i)}^{\text{ref}})}^{-1}}

\newcommand{\KrefinvF}{{(\widehat{\MB{K}}{}_{(i)}^{\text{ref}})}^{-1}}

\newcommand{\KrefinvabF}{{(\widehat{\MB{K}}_{(i)\bar{\alp}\bar{\beta}}^{\text{ref}})}^{-1}}

\newcommand{\linsysMat}{\MB{K}}
\newcommand{\RHS}{\MB{b}}

\newcommand{\mandeld}{\ensuremath{d_{\text{m}}}}

\newcommand{\temperaturesymbol}{\ensuremath{w}}
\newcommand{\temperature}{\temperaturesymbol}
\newcommand{\pertemp}{\Tilde{\temperaturesymbol}}
\newcommand{\Dpertemp}{\MB{\pertemp}}

\newcommand{\conductmaterialsymbol}{A}
\newcommand{\conductmaterial}{\mat{\conductmaterialsymbol}}
\newcommand{\conductDmaterial}{\MB{\conductmaterialsymbol}}
\newcommand{\conductrmaterial}{\mat{\conductmaterialsymbol}^\text{ref}}
\newcommand{\conductDrmaterial}{\MB{\conductmaterialsymbol}^\text{ref}}

\newcolumntype{L}{>{\centering\arraybackslash}m{3cm}}

\usepackage{cancel}  

\begin{document}

\begin{frontmatter}


 \cortext[cor1]{Corresponding author}

\title{Optimal FFT-accelerated Finite Element Solver for Homogenization}


\author[inst1]{Martin Ladecký\corref{cor1}}
\ead{martin.ladecky@cvut.cz}
\author[inst2]{Richard J. Leute}
\ead{richard.leute@imtek.uni-freiburg.de}
\author[inst3]{Ali Falsafi}
\ead{ali.falsafi@epfl.ch}
\author[inst1]{Ivana Pultarová}
\ead{ivana.pultarova@cvut.cz}
\author[inst2]{Lars Pastewka}
\ead{lars.pastewka@imtek.uni-freiburg.de}
\author[inst3]{Till Junge}
\ead{till.junge@epfl.ch}
\author[inst1]{Jan Zeman}
\ead{jan.zeman@cvut.cz}

\address[inst1]{Faculty of Civil Engineering, Czech Technical University in Prague, Th\'akurova 7, 166 29 Prague 6, Czech Republic}

   \address[inst2]{Department of Microsystems Engineering, University of Freiburg, Georges-K\"ohler-Allee 103, 79110 Freiburg, Germany}

 \address[inst3]{Department of Mechanical Engineering, \'Ecole Polytechnique F\'ed\'erale de Lausanne, 1015 Lausanne, Switzerland}

\begin{abstract}
We propose a matrix-free finite element (FE) homogenization scheme that is considerably more efficient than generic FE implementations. The efficiency of our scheme follows from a preconditioned well-scaled reformulation allowing for the use of the conjugate gradient or similar iterative solvers. The geometrically-optimal preconditioner --- a discretized Green’s function of a periodic homogeneous reference problem --- has a block-diagonal structure in the Fourier space which permits its efficient inversion using the fast Fourier transform (FFT) techniques for generic regular meshes. This implies that the scheme scales as $\mathcal{O}(n \log(n))$ like  FFT, rendering it equivalent to spectral solvers in terms of computational efficiency. However, in contrast to classical spectral solvers, the proposed scheme works with FE shape functions with local supports and is free of the Fourier ringing phenomenon. We showcase that the scheme achieves the number of iterations that are almost independent of spatial discretisation and scales mildly with the phase contrast. Additionally, we discuss the equivalence between our displacement-based scheme and the recently proposed strain-based homogenization technique with finite-element projection.

\end{abstract}


\begin{highlights}
\item We propose a specialized preconditioned FE homogenization scheme that is significantly more efficient than generic FE implementations.
\item Our formulation offers computational efficiency equivalent to spectral homogenization solvers.
\item Our formulation is devoid of ringing phenomena intrinsic to spectral methods.
\item All results are reproducible in an open-source project $\mu$Spectre available at \url{https://muspectre.gitlab.io/muspectre}.
\end{highlights}


\end{frontmatter}



\section{Introduction}
Complex macroscopic phenomena such as  plastic yielding or damage in materials are governed by the nonlinear behavior of materials at meso-, micro-, or nanoscales. This intrinsic multiscale aspect of materials behavior creates the demand for the development of specialized scale-bridging techniques \cite{LLorca_2011,MATOUS2017192,Fish2021}.
We focus here on an image-based homogenization technique~\cite{Terada1997} that combines the characterization of materials microstructures by high-resolution images (originating, e.g., from micro-computed tomography~\cite{Maire2014} or geometry-based models~\cite{SONON20211}) and a numerical solution of an underlying partial differential equation (PDE) with coefficient defined on a regular grid and typically involving periodic boundary conditions.


The solution of such PDEs discretized with the conventional finite element (FE) then becomes challenging even in
the simplest scalar elliptic case, because it results in a system of equations with millions to billions of unknowns~\cite[Section 7.6]{johnson_1995}.
In this regard, matrix-free iterative solvers are clearly preferential to direct solvers because of their lower memory footprint and speed, with the conjugate gradient (CG) method~\cite{Hestenes1952MethodsOC} being the most prominent candidate. However, the convergence behavior of the CG method depends on the spectral properties of the linear system matrix and deteriorates with decreasing FE mesh size~\cite[Section 7.7]{johnson_1995}.



More than two decades ago, Moulinec and Suquet in their foundational works \cite{moulinec_fast_1994,moulinec_numerical_1998} proposed a method that resolved these issues. According to its original interpretation, the method employed fixed-point  iterations involving convolution with the Green's function of an auxiliary homogeneous problem with data and unknowns defined directly on the input grid.
The method is suitable for high resolution homogenization problems thanks to the efficient implementation of the convolution step using the fast Fourier transform (FFT) algorithm \cite{golub2013matrix} and mesh-size independent number of iterations.

These features attracted great interest in the community of computational mechanics of materials, as documented in two recent surveys by \citet{Schneider2021} and \citet{Lucarini_2021}. In what follows, we outline the developments most relevant to our work and refer an interested reader to \cite{Schneider2021,Lucarini_2021} for the full story of FFT-based methods.


\emph{Conjugate gradient solvers.} As reported independently by \citet{Brisard2010FFT} and \citet{ZeVoNoMa2010AFFTH}, the original spectral scheme \cite{moulinec_fast_1994,moulinec_numerical_1998} can be further accelerated  when replacing the fixed-point algorithm with the CG method. Later on, these computational observations were justified by \citet{BRISARD2012197}, who showed that the computational scheme of \citet{Brisard2010FFT} follows from the Ritz discretisation of the Hashin-Shtrikman variational principles and by \citet{VoZeMa2014FFTH}, who showed that the computational scheme of \citet{ZeVoNoMa2010AFFTH} follows from the Fourier-Galerkin discretisation of the underlying PDE. These results directly extend to nonlinear problems linearized by the Newton's method, as first reported by \citet{Gelebart2013} and \citet{Kabel2014LargeDef} for the Green's function framework and by \citet{ZeGeVoPeGe2017} and \citet{DEGEUS2017412} for the Fourier-Galerkin framework. 

 
\emph{Oscillations.}  Because the stress or strain fields may exhibit discontinuities at interphases between different material phases, discretizing the problem by Fourier trigonometric polynomials results in spurious numerical oscillations (also referred to as Fourier ringing artifacts in  Section 2.5 of \cite{Schneider2021}) that pollute the approximate results. 
To reduce these oscillations, \citet{Kasbohm2006} smoothed the material data and \citet{SHANTHRAJ201531} filtered out high Fourier frequencies from the solution fields. A different approach was used by \citet{Willot2014}, who considered a modified Green's function obtained from a finite difference discretisation.
 \citet{Schneider2016} extended this approach by proposing a staggered grid finite difference approximation to the underlying PDE, with a follow-up study~\cite{schneider_fft-based_2017} on FE discretisation employing linear hexahedral elements. A related approach building on bi/trilinear FE basis functions instead of the Fourier basis was proposed by 
 \citet{Leuschner2018}. Most recently, \citet{LeuteR2021} developed a compatibility projection-based method in the spirit of Refs.~\cite{ZeGeVoPeGe2017,DEGEUS2017412} while considering several finite difference- and finite element-based discretisation stencils. Further discussion on mitigating the oscillation phenomena can be found in a dedicated comparative study of  \citet{Ma2021} or in Section 2.5 and 2.6 of \citet{Schneider2021}.


\emph{Our work.}
We focus on an alternative FFT-accelerated, oscillation-free computational homogenization scheme based purely on FE discretisation that scales quasilinearly with the mesh size. We consider a nonlinear small-strain elasticity micromechanical problem discretized on a regular periodic grid with FE method in Sections~\ref{sec:small_strain_elasticity} and linearize it with the Newton's method in Sections~\ref{sec:FEM}. Note that the localized support of the FE basis functions directly resolves the oscillation issue, see e.g. \cite{LeuteR2021}. Thus no additional artificial adjustments of the data or the solution are needed.


 In Section~\ref{sec:preconditioning}, we overcome the main drawback of the FE discretisation --- deteriorating conditioning of a linear system with the increasing size of the discretisation grid --- using a suitable preconditioner. Similarly to \cite{schneider_fft-based_2017, Leuschner2018}, we construct the preconditioner from a stiffness matrix of a reference problem with generally anisotropic spatially uniform material data discretized on the same regular grid as the original problem. Using classical results, see e.g.~\cite[Section 5.1.2]{Axelsson2009}, we can guarantee that the condition number of the preconditioned linear system becomes almost independent on the mesh size.  Moreover, employing local ratios of the problem material data and the reference problem material data, we can localize all individual eigenvalues \cite{PultarovaNLAA,Gergelits_2019,Ladecky2020}. This may help to better predict the convergence of the CG method,  see e.g.~\cite[Section 2]{Gergelits_2019}.
 Therefore, the iterative CG solver is an optimal choice for the solution of problems with highly resolved microstructures. The application of the preconditioner is presented in detail in Section~\ref{sec:implementation}, with emphasis on reducing its computational complexity using the FFT algorithm~\cite{CooleyTukey}.



We demonstrate the main features of the proposed algorithm by examples collected in  Section~\ref{sec:experiments} that covers linear thermal conduction (with the necessary adjustments outlined in \ref{sec:appendix_termal}), linear small-strain elasticity, and nonlinear finite strain elasto-plasticity.  Section~\ref{sec:comparison} is devoted to a comparison of our scheme with related developments by \citet{schneider_fft-based_2017} and  \citet{Leuschner2018}, and Section~\ref{sec:conclusion} concludes our work.

\emph{Notation.}
 We denote $d$-dimensional vectors and matrices by boldface letters: $\vek{a}=(a_\alp)_{\alp=1}^{d} \in \D{R}^{d}$ or $ \vek{A}=(A_{\alp \beta})_{\alp,\beta=1}^{d} \in \D{R}^{d \times d} $. Matrix-matrix and matrix-vector multiplications are denoted as $\vek{C}=\vek{B}\vek{A}$ and $\vek{c}=\vek{B}\vek{a}$. 
Vectors and matrices arising from the discretisation will be denoted by $\MB{a} $ and $\MB{A}$, to highlight their special structure. The $(I)$-th component of $\MB{a}$ will be denoted as $\MB{a}[I]$ and $(I,J)$-th component of $\MB{A}$ will be denoted as $\MB{A}[I,J]$.  We consider a general $d$-dimensional setting throughout the paper. However, for the sake of readability, we use $d=2$ in the expanded form of matrices, such as in equation~\eqref{eq:mandel_grad}.

 \section{Nonlinear Small-Strain Elasticity}\label{sec:small_strain_elasticity}
\begin{figure}
\centering
 \includegraphics[width=0.3\paperwidth]{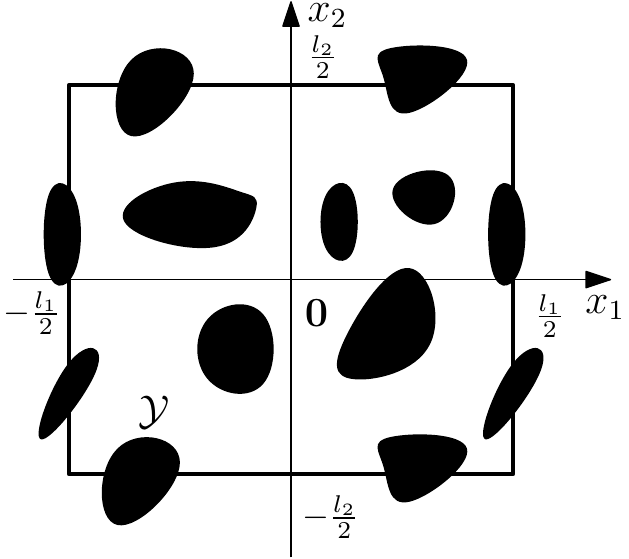}
  \caption{A rectangular two-dimensional cell $\puc= \left[-\dfrac{l_{1}}{2},\dfrac{l_{1}}{2} \right] \times \left[-\dfrac{l_{2}}{2},\dfrac{l_{2}}{2} \right]$ with outlined periodic microstructure.}
  \label{fig:puc}
\end{figure}

We consider a $d$-dimensional rectangular periodic cell $\puc=\prod_{\alp=1}^d {\left[-\frac{l_{\alp}}{2},\frac{l_{\alp}}{2} \right]}$, of volume~$| \puc | =\prod_{\alp=1}^d l_{\alp}$, to be a~representative volume element, i.e.~a~typical material microstructure; see Fig.~\ref{fig:puc} for an illustration. The symmetries of small-strain elasticity allow us to employ the Mandel notation and reduce the dimension of the second-order strain tensor $\symgrad \disp=\frac{1}{2}(\grad{\disp}+\grad{\disp}^{\Transpose}): \puc \rightarrow\D{R}_{\text{sym}}^{d\times d} $ to a vector $\symgradM\disp : \puc \rightarrow \D{R}^{\mandeld}$, where $\symgradM$ is the symmetrized gradient operator  such that, for $d=2$,
\begin{align}\label{eq:mandel_grad}
   \symgradM\disp=
    \begin{pmatrix}
       ( \symgrad \disp)_{11}\\
       ( \symgrad \disp)_{22}\\
       \sqrt{2}( \symgrad \disp)_{12}
    \end{pmatrix}
    =
    \begin{pmatrix}
        \frac{\partial}{\partial x_1}&0\\
        0&\frac{\partial}{\partial x_2}\\
        \frac{\sqrt{2}}{2}\frac{\partial}{\partial x_2}&
        \frac{\sqrt{2}}{2}\frac{\partial}{\partial x_1}
    \end{pmatrix}
    \begin{pmatrix}
        \dispcomponent_1\\
        \dispcomponent_2
    \end{pmatrix}.
\end{align}
Similarly, a fourth-order tensor $ \D{C}: \puc \rightarrow\D{R}_{\text{sym}}^{d\times d\times d\times d} $ is represented with a matrix $\material : \puc \rightarrow \D{R}^{\mandeld \times \mandeld}$,
\begin{align*}
     \material=
        \begin{pmatrix}
        \D{C}_{1111}&\D{C}_{1122}& \sqrt{2} \D{C}_{1112}\\
        \D{C}_{2211}&\D{C}_{2222}& \sqrt{2} \D{C}_{2212}\\
        \sqrt{2}\D{C}_{1211}&\sqrt{2}\D{C}_{1222}& 2        \D{C}_{1212}
    \end{pmatrix},
\end{align*}
 where the number of components of the symmetrized gradient in the Mandel notation is $\mandeld=\cfrac{(d+1)d}{2}$, and indices $\alp_\text{m}, \beta_\text{m}, \gamma_\text{m}\in \left\lbrace  1, \dots ,\mandeld \right\rbrace$.

 In the small-strain micromechanical problem, we split the overall strain $\strain : \puc \rightarrow \D{R}^{\mandeld}$ into an~average strain $\macrostrain=\frac{1}{| \puc |}\int_{\puc} \strain(\vek{x}) \intd{\vek{x}}  \in \D{R}^{\mandeld}$ and a periodically fluctuating field $\symgradM \perdisp : \puc \rightarrow \D{R}^{\mandeld}$,
\begin{align*}
    \strain(\vek{x})
    =
    \macrostrain+\symgradM \perdisp(\vek{x}) \quad \text{for all }  \vek{x} \in \puc.
\end{align*}
Here,~$\symgradM \perdisp $ denotes the symmetrized gradient in the Mandel notation, and the fluctuating displacement field $\perdisp$~belongs to the space of admissible functions $
\testspace=\left\lbrace  \perdisptest :\puc\rightarrow\D{R}^{d},\, \perdisptest \text{ is } \puc \text{-periodic} \right\rbrace$.
The governing equations for $\symgradM \perdisp$ are the mechanical equilibrium conditions
 \begin{align*}
     -\symgradM^{\Transpose}  \stress(\vek{x},\macrostrain+\symgradM \perdisp(\vek{x}) ,\vek{g}(\vek{x}))=\vek{0} \quad \text{for all }  \vek{x} \in \puc,
 \end{align*}
in which $\stress:\puc \times \D{R}^{\mandeld}\times \D{R}^{g}\rightarrow \D{R}^{\mandeld}$ is the stress field and $\vek{g}:\puc\rightarrow\D{R}^{g}$ designates the vector of internal parameters. The equilibrium equations are converted to the weak form
\begin{align}\label{eq:weak_form}
    \int_{\puc}
   \symgradM \perdisptest(\vek{x})^{\Transpose}
    \stress(\vek{x},\macrostrain+\symgradM \perdisp(\vek{x}),\vek{g}(\vek{x}) )
    \intd{\vek{x}}
    = 0  \quad \quad \text{for all }  \perdisptest \in \testspace,
\end{align}
where $\perdisptest$ is the test displacement field. The weak form \eqref{eq:weak_form} serves as a starting point for the FE method.
 
 \section{Finite Element Discretisation}\label{sec:FEM}
\begin{figure*}[htbp]
  \centering
 \includegraphics[width=1.0\textwidth]{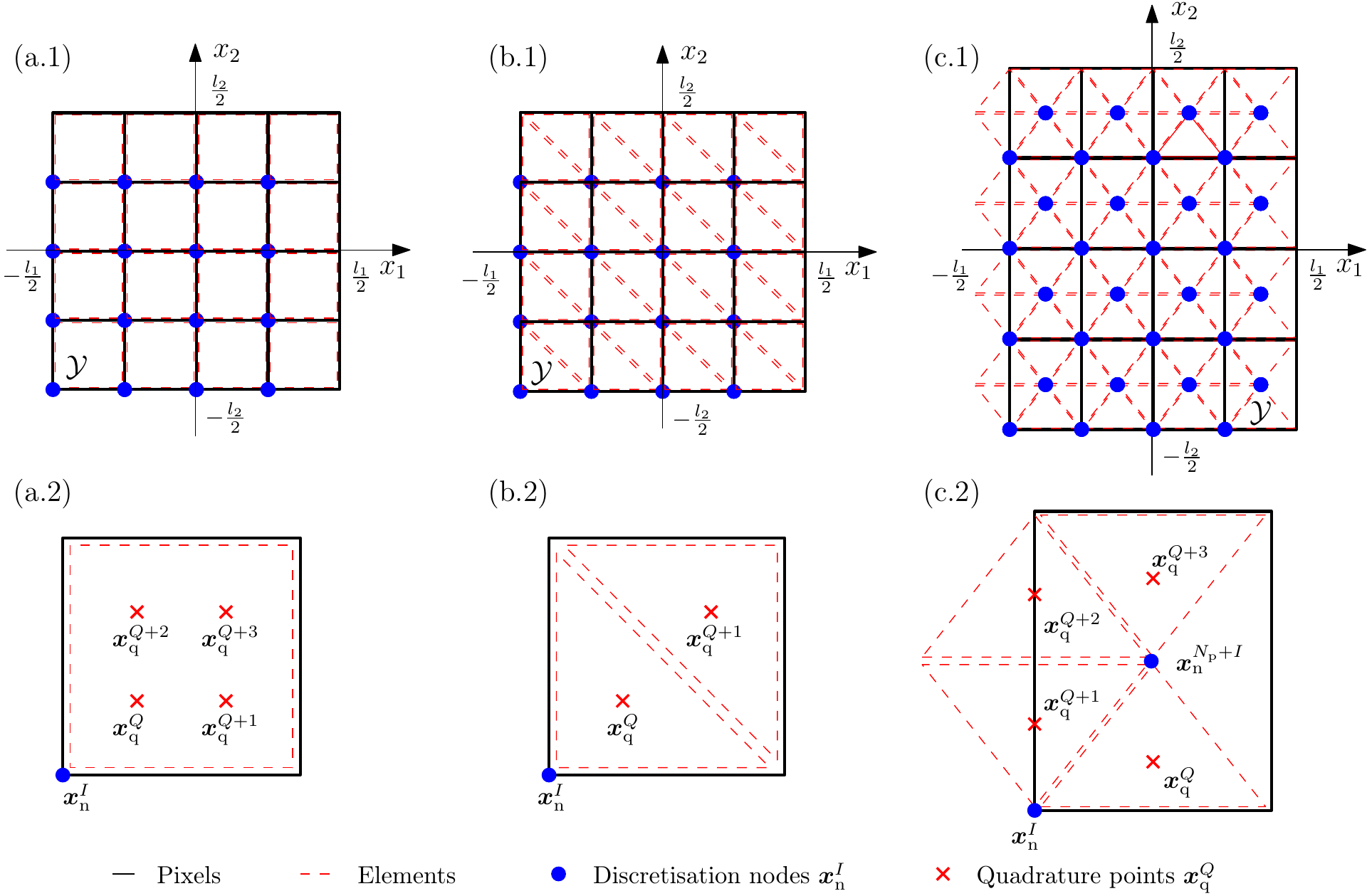}
  \caption{Example of regular periodic FE grids with associated discretisation stencils for a two-dimensional~cell~$\puc$. All grids consists of $16$ pixels ($\Np=16$).
  The row \textbf{(1)} shows:
  \textbf{(a.1)} grid with $16$ discretisation nodes ($\NI=16$) and quadrature points  ($\NQ=64$),
  \textbf{(b.1)} grid with $16$ discretisation nodes ($\NI=16$) and $32$ quadrature points  ($\NQ=32$), \textbf{(c.1)} grid with $32$ discretisation nodes  ($\NI=32$) and $64$ quadrature points  ($\NQ=64$).
  The row \textbf{(2)} shows:
  \textbf{(a.2)} one-node stencil ($\Nn=1$) with one bilinear rectangular element and four quadrature points with the quadrature weights $w^Q=\frac{1}{4}V_{\text{p}}$, \textbf{(b.2)} one-node stencil ($\Nn=1$) with two linear triangular elements and two quadrature points with the quadrature weights $w^Q=\frac{1}{2}V_{\text{p}}$, \textbf{(c.2)} two-node stencil ($\Nn=2$) with four linear triangular elements and four quadrature points with the quadrature weights $w^Q=\frac{1}{4}V_{\text{p}}$,  Here, $V_{\text{p}}$ denotes pixel volume, such that $V_{\text{p}}\Np=| \puc |$.
}
  \label{fig:tiling}
\end{figure*}

For the discretisation of the weak form \eqref{eq:weak_form}, we use a uniform mesh and conforming FE basis functions.
In our setting, the discretisation mesh does not necessarily follow the regular pixel/voxel structure, but can correspond to a space-filling pattern of finite elements; see the first row in Fig.~\ref{fig:tiling}.
 The discretisation mesh is generated by a periodic repetition of a discretisation stencil in the cell $\puc$; see the second row in Fig.~\ref{fig:tiling}. Such flexibility in discretisation is useful, e.g., for material models that exhibit sensitivity to mesh-grid anisotropy.

Strain and stress fields are evaluated at quadrature points $\vek{x}_{\mathrm{q}}^Q$,~$ Q\in \{1, 2, \dots , \NQ\, \}$, cf.~Fig.~\ref{fig:tiling}, and
 the displacement fields are sampled at discretisation nodes $\vek{x}_{\mathrm{n}}^{I} $,~$ I\in \{1, 2, \dots , \NI \, \}$. The number of discretisation nodes $\NI=\Np \Nn$ is given by the number of pixel/voxel-associated discretisation stencils $\Np$ and the number of nodes per stencil $\Nn$, as explained in Fig.~\ref{fig:tiling}. The number of degrees of freedom per stencil is thus $d\Nn$ and the total number of degrees of freedom per domain is~$d\NI$.
 
%
Following the standard FE theory, $\perdisptest$ and $\vek{\Tilde{u}}$ are approximated by continuous element-wise polynomials $\mathcal{P}_k$ of the degree $k$; their symmetrized gradients $\symgradM \perdisptest$ and $\symgradM\perdisp$ then become element-wise polynomials of the degree up to~$k$.
Furthermore, the integral \eqref{eq:weak_form} can be approximated with a suitable quadrature rule,
\begin{equation}\label{eq:global_sum}
\begin{split}
        \int_{\puc}	 \symgradM \perdisptest(\vek{x})^{\Transpose}
    &\stress(\vek{x},\macrostrain+\symgradM \perdisp (\vek{x}) ,\vek{g}(\vek{x}))
    \intd{\vek{x}}\\
        \approx 
        \sum_{Q=1}^{\NQ} 
         \symgradM   \perdisptest (\vek{x}_{\mathrm{q}}^Q)^{\Transpose}
        &\stress(\vek{x}_{\mathrm{q}}^Q,\macrostrain+\symgradM  \perdisp(\vek{x}_{\mathrm{q}}^Q),\vek{g}(\vek{x}_{\mathrm{q}}^Q) ) \,
        w^Q,
\end{split}
\end{equation}
 where the positions of the quadrature points $\vek{x}_{\mathrm{q}}^Q$ and the quadrature weights $w^Q$ depend on the choice of the quadrature rule; recall~Fig.~\ref{fig:tiling}.
 
 Every component $ \perdispcomponent_\alp $ of the unknown vector $\perdisp$ is approximated by a linear combination
    \begin{align*}
        \perdispcomponent_\alp(\vek{x}) \approx \perdispcomponent_{\alp}^{N}(\vek{x}) 
        = \sum_{I=1}^{\NI}
        \perdispcomponent_{\alp}^{\N}(\vek{x}_{\mathrm{n}}^{I})\phi^{I}(\vek{x}) 
        \quad \text{for all }  \vek{x} \in \puc,
    \end{align*}
    where the coefficients $\perdispcomponent_{\alp}^{\N}(\vek{x}_{\mathrm{n}}^{I})$ are the nodal values of $\perdispcomponent_{\alp}^{\N}$ at discretisation nodes $\vek{x}_{\mathrm{n}}^{I}$ and $\phi^{I}$ are FE basis functions.
    A partial derivative of this approximation
    \begin{align*}
          \frac{\partial \perdispcomponent_{\alp}^{\N}(\vek{x}) }{\partial x_\beta}  
          =
          \sum_{I=1}^{\NI}
            \perdispcomponent_{\alp}^{\N}(\vek{x}_{\mathrm{n}}^I) 
            \frac{\partial \phi^{I}(\vek{x}) }{\partial x_\beta}  
            \quad \text{for all }  \vek{x} \in \puc,
    \end{align*}
    evaluated in the quadrature points is given by
    \begin{align*}
       \frac{\partial \perdispcomponent_{\alp}^{\N}(\vek{x}_{\mathrm{q}}^Q) }{\partial x_\beta }  =\sum_{I=1}^{\NI}
        \perdispcomponent_{\alp}^{\N}(\vek{x}_{\mathrm{n}}^I)
        \frac{\partial \phi^{I}(\vek{x}_{\mathrm{q}}^Q)}{\partial x_\beta}
        \quad \text{for } Q =1, \dots , \NQ
        .
 \end{align*}
  Therefore, if we store the nodal values of displacement $ \perdisp(\vek{x}_{\mathrm{n}}^I)$ into a vector $\MB{ \Dperdisp} \in \D{R}^{d\NI} $, the gradient vector $\MB{\symgradM \Dperdisp}\in \D{R}^{\mandeld\NQ} $ at all quadrature points is given with 
    \begin{align}\label{eq:Du_mat}
\symgradM \Dperdisp
    =
    \Dmat\MB{\Dperdisp}=
    {\begin{bmatrix}
   \Dmat_{1} & \MB{0}    \\ 
   \MB{0}   & \Dmat_{2}  \\
  \frac{\sqrt{2}}{2} \Dmat_{2} &  \frac{\sqrt{2}}{2}\Dmat_{1}    
   \end{bmatrix}}
    {\begin{bmatrix}
   \MB{\Dperdisp}_{1}  \\
   \MB{\Dperdisp}_{2} 
   \end{bmatrix}},
    \end{align}
where the matrix $\Dmat\in\D{R}^{\mandeld\NQ \times d\NI}$ consists of sub-matrices of the partial derivatives
 \begin{align}\label{eq:grad_componets}
    \MB{D}_{\beta}[Q, I]=\frac{\partial \phi^{I}(\vek{x}_{\mathrm{q}}^Q
     )}{\partial x_\beta}
     \quad \text{for } Q =1, \dots , \NQ
     \text{ and }
      I =1, \dots , \NI
     ,
 \end{align}
 and $ \MB{\Dperdisp}_{\alpha}$ stores values of the displacement in the direction ${\alpha} $.
Due to the local supports of the basis functions $\phi^{I}$, these sub-matrices exhibit significant sparsity, e.g., for the element-wise linear approximation, shown in the middle of Fig.~\ref{fig:tiling}, each row of $\Dmat_{\beta}$ contains only two nonzero entries. Since both the interpolating and quadrature points are periodically distributed in~$\puc$, the matrix $\Dmat_{\beta}$ has a block circulant structure.

Now, the discretized weak form \eqref{eq:weak_form} using quadrature \eqref{eq:global_sum} can be rewritten in the matrix notation as 
         \begin{align}\label{eq:disc_wf}
          \MB{\testf}^{\Transpose}\Dmat^{\Transpose}
           \MB{\qweights}
          \stress(\MB{e}+
         \Dmat\MB{\Dperdisp},\MB{g})
         =
         0
         \quad \text{for all } \MB{\testf} \in \D{R}^{d\NI}  ,
    \end{align}
where $\MB{\testf}$ stores the nodal values of test displacements, $\MB{e} \in \D{R}^{\mandeld\NQ} $ stands for the discretized average strain, $  \stress: \D{R}^{\mandeld\NQ} \times  \D{R}^{g\NQ} \rightarrow \D{R}^{\mandeld\NQ}$ is a nonlinear map transforming, locally at quadrature points, a vector of discrete strains and internal parameters $\MB{g}\in \D{R}^{g\NQ}$ to
discrete stresses, and the diagonal matrix $\MB{\qweights}\in \D{R}^{\mandeld\NQ\times \mandeld\NQ} $
    \begin{align}\label{eq:W_mat}
        \MB{\qweights}
        =
        {\begin{bmatrix}
           \MB{\qweights_\text{m}} & \MB{0} & \MB{0}   \\ 
           \MB{0}   & \MB{\qweights_\text{m}} & \MB{0}          \\ 
           \MB{0}   & \MB{0}    & \MB{\qweights_\text{m}}
       \end{bmatrix}}
   \end{align}
consists of $\mandeld$ identical diagonal matrices $\MB{\qweights_\text{m}}\in \D{R}^{\NQ\times \NQ}$ storing quadrature weights, $\MB{\qweights_\text{m}}[Q,Q]=w^Q$.

As the vector $\MB{\testf}$ is arbitrary, discretized weak form \eqref{eq:disc_wf} is equivalent to the system of  discrete nonlinear  equilibrium conditions
    \begin{align}\label{eq:non_lin_system}
          \Dmat^{\Transpose}
           \MB{\qweights}
          \stress(\MB{e}+
         \Dmat\MB{\Dperdisp},\MB{g})
         =
         \MB{0}.
    \end{align}
\subsection{Linearisation}\label{sec:linearisation}
We employ the Newton's method to solve the nonlinear system \eqref{eq:non_lin_system} iteratively. To this purpose, the $(i+1)$-th approximation of the nodal displacement $\MB{\Dperdisp}_{(i+1)}\in \D{R}^{\NI}$ is given by the previous approximation $\MB{\Dperdisp}_{(i)}\in \D{R}^{\NI}$ adjusted by a finite displacement increment~$\delta\MB{\Dperdisp}_{(i+1)}\in \D{R}^{\NI}$,
\begin{align*}
    \MB{\Dperdisp}_{(i+1)}
     =
    \MB{\Dperdisp}_{(i)}+\delta\MB{\Dperdisp}_{(i+1)},
\end{align*}
with an initial approximation $\MB{\Dperdisp}_{(0)}\in \D{R}^{\NI}$. The displacement increment $\delta\MB{\Dperdisp}_{(i+1)}$ follows from the solution of the linear system
\begin{align}\label{eq:lin_system}
\underbrace{
    \Dmat^{\Transpose}
     \MB{\qweights}
    \MB{C}_{(i)} \Dmat
     }_{\linsysMat_{(i)}}
    \delta\MB{\Dperdisp}_{(i+1)}
    =
     \underbrace{
    -
    \Dmat^{\Transpose}
     \MB{\qweights}
    \stress(\MB{e}+
    \Dmat \MB{\Dperdisp}_{(i)}, \MB{\vek{g}}_{(i)})
     }_{\RHS_{(i)}},
\end{align}
where the discrete constitutive tangent matrix           $
          \Dmaterial_{(i)}
         =
        \dfrac{\partial \stress}{\partial \strain}(\MB{e}+
         \Dmat \MB{\Dperdisp}_{(i)}, \MB{\vek{g}}_{(i)}) \in  \D{R}^{\mandeld\NQ\times \mandeld\NQ}, 
         $
         \begin{align*}
         \Dmaterial_{(i)} =
  {\begin{bmatrix}
   \Dmaterial_{(i) 11 }& \Dmaterial_{(i)12}  & \Dmaterial_{(i)13}  \\
   \Dmaterial_{(i) 21} & \Dmaterial_{(i)22}  & \Dmaterial_{(i)23}  \\
   \Dmaterial_{(i) 31} & \Dmaterial_{(i)32}  & \Dmaterial_{(i)33} 
   \end{bmatrix}},
    \end{align*}
is obtained from the constitutive tangent 
 $\material_{(i)}(\vek{x}) =\dfrac{\partial \stress}{\partial \strain}(\vek{x},\macrostrain+\symgradM \perdisp_{(i)} (\vek{x}),\vek{g}_{(i)}(\vek{x}))$, evaluated at quadrature points.
Therefore, the sub-matrices $\Dmaterial_{(i)\alpha_\text{m}\beta_\text{m}}\in  \D{R}^{\NQ\times \NQ}$ are diagonal with entries $\Dmaterial_{(i)\alpha_\text{m}\beta_\text{m}}[Q,Q]=C_{(i)\alpha_\text{m}\beta_\text{m}}(\vek{x}_{\mathrm{q}}^Q)$.
Traditionally, $\linsysMat_{(i)}\in  \D{R}^{d\NI\times d\NI}$  denotes the matrix of the linear system \eqref{eq:lin_system}, and $\RHS_{(i)}\in  \D{R}^{d\NI}$ stands for the right-hand side of \eqref{eq:lin_system}.

 \section{Preconditioning}\label{sec:preconditioning}
 Recall that we focus on micromechanical problems with a finely described microstructure that involves a large number of degrees of freedom $d\NI$. We aim to use a memory-efficient matrix-free iterative method to find the solution of the linear system \eqref{eq:lin_system}.
 The system matrix $\linsysMat_{(i)}$ is symmetric and positive definite for the symmetric constitutive tangent $\Dmaterial_{(i)}$, which renders the CG method as the method of choice, when combined with an appropriate preconditioner.
This section discusses how to construct such a preconditioner in an optimal manner.
 
 
\subsection{Reference Material-Based Preconditioner}
The idea of preconditioning, see, e.g.,~\cite[Section~10.3]{golub2013matrix}
or~\cite[Chapters~9 and~10]{Saad2003}, is based on assumptions that the matrix of the preconditioned linear system
\begin{align}\label{eq:prec_lin_system_general}
\MB{M}_{(i)}^{-1}\linsysMat_{(i)} \delta\MB{\Dperdisp}_{(i+1)}=\MB{M}_{(i)}^{-1}\RHS_{(i)},
\end{align}
has more favourable spectral properties 
than the original system $\linsysMat_{(i)} \delta\MB{\Dperdisp}_{(i+1)}=\RHS_{(i)}$. At the same time, the preconditioning matrix $ \precMatrix_{(i)}\in  \D{R}^{d\NI\times d\NI}$ should be relatively easy to invert, such that the faster convergence of the iterative method compensates the computational overhead of the preconditioning. 
Please note that system matrix $\MB{M}_{(i)}^{-1}\linsysMat_{(i)}$ is no longer symmetric. However, for symmetric $\MB{M}_{(i)}$ and $\linsysMat_{(i)}$, system~\eqref{eq:prec_lin_system_general} is equivalent with the system preconditioned in the symmetric form $\MB{M}_{(i)}^{-1/2}\linsysMat_{(i)}\MB{M}_{(i)}^{-1/2} \delta\MB{z}_{(i+1)}=\MB{M}_{(i)}^{-1/2}\RHS_{(i)}$, where $\delta\MB{z}_{(i+1)} =\MB{M}_{(i)}^{1/2}\delta\MB{\Dperdisp}_{(i+1)}$. The latter form is in fact solved by the PCG method; see~\cite[Section~9.2.1]{Saad2003} for more details. Nonetheless, we prefer the notation with the left preconditioning \eqref{eq:prec_lin_system_general} for brevity.

Our approach is based on a preconditioner constructed in the same manner as the original matrix of the linear system~$\eqref{eq:lin_system}$,
  \begin{align}\label{eq:Kref}
\MB{M}_{(i)}=\Kref =
  \Dmat^{\Transpose}  \MB{\qweights}\Drmaterial_{(i)}\Dmat\in  \D{R}^{d\NI\times d\NI},
 \end{align}
where the reference constitutive tangent matrix $\Drmaterial_{(i)} \in  \D{R}^{\mandeld\NQ\times \mandeld\NQ}$ corresponds to spatially uniform material data $\rmaterial_{(i)} \in  \D{R}^{\mandeld\times \mandeld}$. 
Finally, substituting \eqref{eq:Kref} into \eqref{eq:prec_lin_system_general} leads to the preconditioned linear system 
\begin{align}\label{eq:prec_lin_system}
\Krefinv \linsysMat_{(i)} \delta\MB{\Dperdisp}_{(i+1)}=\Krefinv\RHS_{(i)},
\end{align}
referred to as the reference material-based preconditioned problem in what follows. Notice that the spectrum of $\Kref $ contains null eigenvalue(s), associated with the rigid body modes, thus instead of the inverse of $\Kref$, we consider its (Moore-Penrose) pseudo-inverse\footnote{For details of Moore-Penrose pseudo-inverse refer to~\cite{golub2013matrix}} but still denote it by $\Krefinv$ for notation simplicity.

 In the following, we advocate this choice of the preconditioner. First, we derive a computationally efficient pseudo-inverse of $\Kref$ and second, we explain how the preconditioning impacts the spectral properties of the matrix of the system~\eqref{eq:prec_lin_system}.

\subsection{Fourier Pseudo-Inversion of $\Kref$}

\begin{figure}
\centering
 \includegraphics[width=1\textwidth]{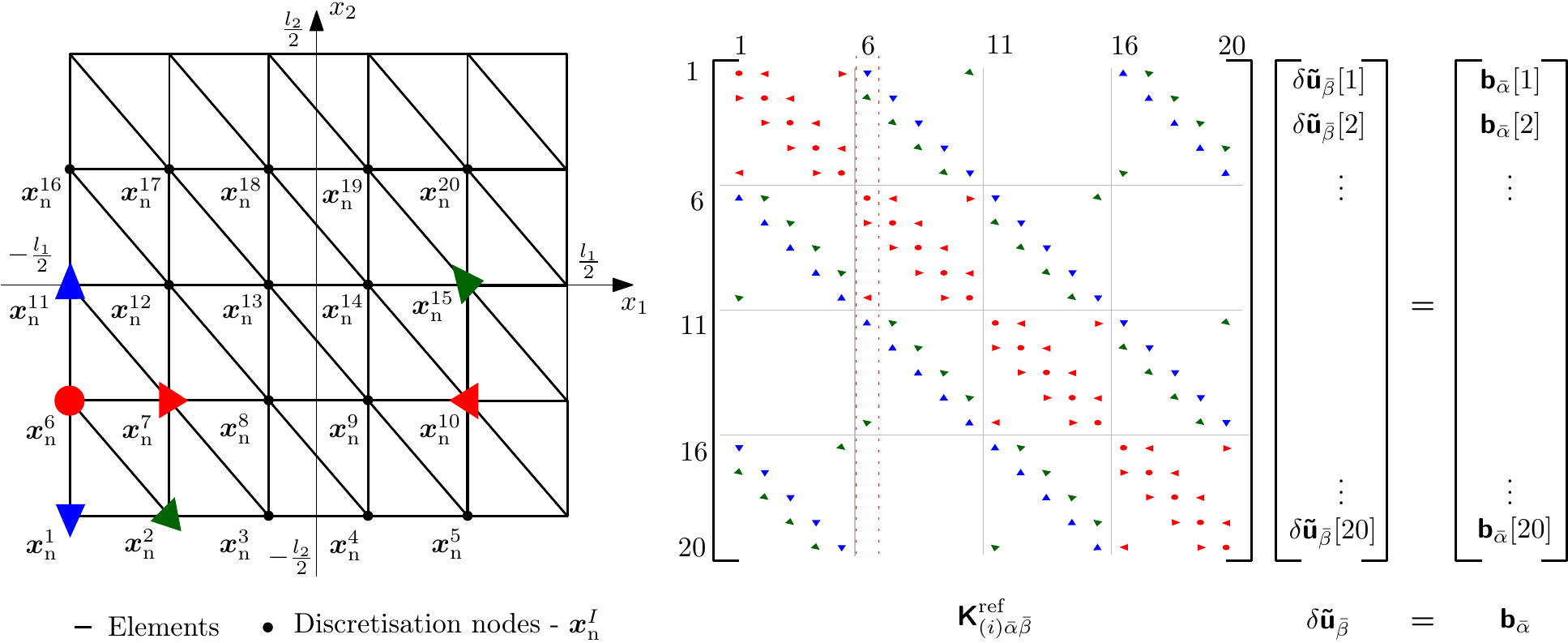}
  \caption{The block-circulant structure of block $\Krefab$ from the preconditioner $\Kref$ for spatially uniform material data $\rmaterial_{(i)}$ and periodic boundary condition. The two-dimensional ($d=2$) discretisation grid  consisting of $20$ pixels ($\Np=20$) with one-node stencil ($\Nn=1$), and $20$ discretisation nodes ($\NI=20$) is shown left. Contributions of unit nodal displacement
$\delta\boldsymbol{\mathsf{\tilde{u}}}_{\bar{\beta}}[I]=1$ to nodal components of right-hand side vector, graphically shown in the node $\boldsymbol{\vek{x}}_{\mathrm{n}}^{6}$,
are given as follows:
(\textcolor{red}{$\bullet$})~ self contribution, contributions
(\textcolor{red}{$\blacktriangleright$})~to the right node,
(\textcolor{red}{$\blacktriangleleft$})~to the left node, 
(\textcolor{OliveGreen}{\protect\rotatebox[origin=c]{-45}{$\blacktriangleleft$}})~to the upper left node,
 (\textcolor{blue}{\protect\rotatebox[origin=c]{-90}{$\blacktriangleleft$}})~to the upper node,
 (\textcolor{blue}{\protect\rotatebox[origin=c]{90}{$\blacktriangleleft$}})~to the bottom node, and
  (\textcolor{OliveGreen}{\protect\rotatebox[origin=c]{135}{$\blacktriangleleft$}})~to the bottom  right node. 
 }
  \label{fig:convolution}
\end{figure}

Regular FE discretisation of the problem with periodic boundary conditions leads to the same stencil for every pixel. 
Thus, for the uniform $\rmaterial_{(i)}$ in the whole $\puc$ (at every quadrature point $\vek{x}_{\mathrm{q}}^Q$), the resulting preconditioning matrix $ \Kref 
  \in \D{R}^{d\Nn\Np\times d\Nn\Np}$, 
  \begin{align}\label{eq:M}
 \Kref 
  =
  \begin{bmatrix}
   \Krefcomp{1}{1} & \Krefcomp{1}{2} \\
     \Krefcomp{2}{1} & \Krefcomp{2}{2} 
   \end{bmatrix}\in \D{R}^{2\Np\times 2\Np},  \quad  (\text{for } \, d\Nn=2)
 \end{align}
 consists of $(d\Nn)^2$ block-circulant blocks $\Krefab \in \D{R}^{\Np\times \Np}$, where $\bar{\alp},\bar{\beta} \in \{ 1, \dots, d \Nn \} $. All row vectors of a block-circulant block $\Krefab$ contain the same information and each row is block-periodically shifted with respect to the preceding one. This directly reflects the periodically repeated discretisation pattern; recall Fig.~\ref{fig:tiling}, and that the action of $\Krefab$ is a discrete convolution of the displacement $\delta\boldsymbol{\mathsf{\tilde{u}}}_{\bar{\beta}}$ with the discretisation kernel, as schematically shown
 in~Fig.~\ref{fig:convolution}.
 Note that in the one-dimensional ($d=1$) case with one node per interval ($\Nn=1$), $\Kref$ has only one circulant block, $\Kref=\Krefcomp{1}{1}$. The block structure of $\Kref $ appears whenever more than one type of degree of freedom is involved, i.e., $d>1$, or $\Nn > 1$.
 
To make the inversion of $\Kref$ efficient, let us define the discrete $d$-dimensional Fourier transform matrix $\MB{F} \in \D{R}^{\Np \times \Np}$ such that  $\MB{F}^{\HTran}=\MB{F}^{-1}$, where $\MB{F}^{\HTran}$ is the conjugate transpose of~$\MB{F}$. Then the Fourier counterpart 
\begin{align*}
\hatKrefab
=
\MB{F}\Krefab\MB{F}^{\HTran}
\end{align*} 
to any block-circulant $\Krefab$ is diagonal, and has the same spectrum (eigenvalues) as $\Krefab$.
Therefore, $\KrefF $ is block-diagonal and cheaply (pseudo) invertible
  \begin{align}\label{eq:M_inverse}
   \Krefinv
   =
 \MB{F}_d^{\HTran}
\KrefinvF  \MB{F}_d
   =
     \begin{bmatrix}
   \MB{F}^{\HTran} & \MB{0} \\
     \MB{0}&\MB{F}^{\HTran}
   \end{bmatrix}
   {
  \begin{bmatrix}
   \KrefFcomp{1}{1}& \KrefFcomp{1}{2}    \\
   \KrefFcomp{2}{1} &
   \KrefFcomp{2}{2}
   \end{bmatrix}^{-1} 
        \begin{bmatrix}
   \MB{F} & \MB{0}    \\
     \MB{0}&\MB{F} 
   \end{bmatrix}
},
 \end{align}
where $\MB{F}_d=I_{d \Nn} \otimes \MB{F} $ and $ I_{d \Nn} \in \D{R}^{d \Nn \times d \Nn}$ is the identity matrix. The expanded form in \eqref{eq:M_inverse} apply for $d\Nn=2$.

Finally, inserting \eqref{eq:M_inverse} as the preconditioner in \eqref{eq:prec_lin_system} leads to
    \begin{align}\label{eq:prec_lin_system_FFT} 
    \underbrace{
    \MB{F}_d^{\HTran}
    \KrefinvF \MB{F}_d}_{\Krefinv}
         \linsysMat_{(i)}\delta\Dperdisp_{(i+1)}
         =&
         \underbrace{
         \MB{F}_d^{\HTran}\KrefinvF \MB{F}_d
        }_{\Krefinv}
         \MB{\RHS}_{(i)},
\end{align}
which reads in the expanded form as

\begin{equation}\label{eq:prec_scheme}
\begin{split}
        &\underbrace{
        \MB{F}_d^{\HTran}(\MB{F}_d\Dmat^{\Transpose} \MB{\qweights} \MB{\Drmaterial_{(i)}} \Dmat\MB{F}_d^{\HTran} )^{-1}\MB{F}_d
        }_{\Krefinv}
         \underbrace{
         \Dmat^{\Transpose}
          \MB{\qweights}
        \MB{C}_{(i)}
         \Dmat
          }_{\linsysMat_{(i)}}
         \delta\Dperdisp_{(i+1)}
         \\&=
         -
          \underbrace{
         \MB{F}_d^{\HTran}(\MB{F}_d\Dmat^{\Transpose}  \MB{\qweights}\MB{\Drmaterial_{(i)}} \Dmat\MB{F}_d^{\HTran} )^{-1}\MB{F}_d
         }_{\Krefinv}
          \underbrace{
         \Dmat^{\Transpose} \MB{\qweights} \stress(\MB{e}+\Dmat \Dperdisp_{(i)}, \MB{\vek{g}}_{(i)})
         }_{-\RHS_{(i)}}.
\end{split}
\end{equation}

\subsection{Spectrum of the Preconditioned Problem}\label{sec:condition_number}
To support the claim that the system matrix of the linear system \eqref{eq:prec_lin_system_FFT} is well conditioned, we rely on the results published recently in~\cite{PultarovaNLAA, Gergelits_2019, Ladecky2020} that provide simple algorithms for obtaining guaranteed two-sided bounds for all individual eigenvalues of the preconditioned operator by using element-by-element estimates. Note that extremal eigenvalue bounds obtained by such an element-by-element algorithm were introduced first in~\cite{Axelsson2009,Eijkhout1991TheRO} and found use,  e.g., in algebraic multilevel methods~\cite{axelsson_1996}. 
Recently, motivated by \citet{nielsen_preconditioning_2009}, \citet{Gergelits_2019} published a new method yielding the bounds to all individual eigenvalues. This allows not only estimating the condition number of the preconditioned system but also to characterize its spectrum, which can provide more specific insights into the convergence of the CG method; see e.g.~\cite[Section 2]{Gergelits_2019} for more details. In \cite{PultarovaNLAA,Ladecky2020} an alternative algorithm is presented that can be applied to a variety of problems and discretisation methods.

Let us recall the approach of \citet{Ladecky2020}. Thanks to the local supports of FE basis functions $\phi^{I}$ it is possible to estimate all eigenvalues of the preconditioned linear system matrix \eqref{eq:prec_lin_system_FFT}. For each $\phi^{I}$, we calculate
\begin{align*}
\lambda_{I}^{\rm L}&=\min_{\vek{x}_{\mathrm{q}}^Q \in \supp\phi^{I}}\lambda_{\rm min}
    \left((\rmaterial_{(i)}(\vek{x}_{\mathrm{q}}^Q))^{-1}{\material_{(i)}}(\vek{x}_{\mathrm{q}}^Q)\right), \quad
  I=1,\dots, \NI,\; \\
   \lambda_{I}^{\rm U}
   &=
   \max_{\vek{x}_{\mathrm{q}}^Q \in \supp\phi^{I}}\lambda_{\rm max}
    \left((\rmaterial_{(i)}(\vek{x}_{\mathrm{q}}^Q))^{-1}{\material_{(i)}}(\vek{x}_{\mathrm{q}}^Q)\right), \quad
  I=1,\dots, \NI,
\end{align*}
where $\supp\phi^{I}$ denotes the support of $\phi^{I}$, and $\lambda_{\rm min},\lambda_{\rm max} $ are the  minimal and  maximal generalized eigenvalues, respectively. For element-wise constant materials $\material_{(i)}$ and $\rmaterial_{(i)}$, any quadrature point $\vek{x}_{\mathrm{q}}^Q$ can be used to evaluate $\lambda_{\rm min}$ and $\lambda_{\rm max}$ on element. Therefore, only one pair $\lambda_{\rm min},\lambda_{\rm max}$ has to be calculated for each element.
Considering every $\lambda_{I}^{\rm L}$ and $\lambda_{I}^{\rm U}$ $d$-times and sorting these two sets into nondecreasing sequences gives the desired lower and upper eigenvalue bounds.

The resulting eigenvalue bounds are therefore independent of the characteristic element diameter $h$, which suggests that the condition number\footnote{Please note that by the condition number $\kappa(\Krefinv \linsysMat_{(i)} )$ we mean the ratio of the largest and the smallest eigenvalues of $\Krefinv \linsysMat_{(i)}$.} $\kappa(\Krefinv \linsysMat_{(i)} )$ of the preconditioned linear system \eqref{eq:prec_lin_system_FFT} will be independent of the problem size. In contrast, $\kappa(\linsysMat_{(i)} )=\mathcal{O}(h^{-2})$ for the unpreconditioned problem, e.g.~\cite[Section~7.7]{johnson_1995}. The ratio between the maximum and minimum eigenvalues of the preconditioned problem \eqref{eq:prec_lin_system_FFT} will increase with an increasing ratio between extreme eigenvalues of $\material_{(i)}$ (so-called material contrast) and decrease as the reference material data $\rmaterial_{(i)}$ approach the material data $\material_{(i)}$ of the problem. Therefore, we can call our preconditioner as optimal, or more precisely, as \emph{geometrically optimal}, which emphasizes that by keeping the discretisation and changing only the data of the preconditioner can lead to the matrix where all eigenvalues are the same, i.e. the condition number is $1$. However, in such a case, the inversion of the preconditioner would become more expensive. The effects of phase contrast and the choice of $\rmaterial_{(i)}$ on the CG performance are further illustrated with examples presented later in Section~\ref{sec:example_small_elasticity}.

 \section{Implementation}\label{sec:implementation}
 The pseudo-algorithm of the incremental Newton-PCG solver for FE discretisation on a regular grid is outlined in Algorithm~\ref{alg:main}.
In the first part, we detail a matrix-free implementation. The second part deals with the assembly of the preconditioner via matrix-free operators and the third part focuses on the efficient pseudo-inversion of the preconditioner.

\subsection{Matrix-free Implementation}
As mentioned in the previous sections, the explicit matrix structure is useful for explanation, but the computations can be performed more efficiently in a matrix-free manner. 

\emph{The Gradient.}
Computational efficiency of our method relies on the fast evaluation of the gradient vector $
  \MB{\symgradM \Dperdisp}
    =
    \Dmat\MB{\Dperdisp}
$. For regular periodic discretisations, the multiplication $\Dmat\MB{\Dperdisp}$ can be implemented as a convolution of $\MB{\Dperdisp}$ with a short kernel, namely the gradient stencil. To emphasize this, we replace matrix notation $\Dmat$ and $\Dmat^{\Transpose}$ with the (matrix-free) operator notation $\Doper:\D{R}^{d\NI} \rightarrow \D{R}^{\mandeld\NQ}$  and $\DoperTrans:\D{R}^{\mandeld\NQ} \rightarrow \D{R}^{d\NI}$, such that
 \begin{align*}
      \Doper \delta\Dperdisp_{(i+1)}=\Dmat \delta\Dperdisp_{(i+1)}, 
      \quad
      \text{and}
      \quad
         \DoperTrans\MB{\qweights}
        \MB{\Dmaterial}_{(i)}
         \Doper
         \delta\Dperdisp_{(i+1)}
         =
      \Dmat^{\Transpose}  \MB{\qweights}
        \MB{\Dmaterial}_{(i)}
         \Dmat
         \delta\Dperdisp_{(i+1)}
         .
\end{align*}
 These operations are equal from the viewpoint of linear algebra, but algorithmically $\Doper$ is of linear $\mathcal{O}(\NI)$ cost.

\emph{The fast Fourier transform.}
In the same manner, the multiplication with the discrete Fourier transform matrix can be replaced with the forward and the inverse fast Fourier transform algorithm  
 \begin{align*}
\FourierTrans \delta\Dperdisp_{(i+1)}= \MB{F} \delta\Dperdisp_{(i+1)}
    \quad
      \text{and}
      \quad
\iFourierTrans \delta\Dperdisp_{(i+1)}=\MB{F}^{H} \delta\Dperdisp_{(i+1)},
\end{align*}
of $\mathcal{O}(\NI\log \NI)$ complexity. 

\emph{Quadrature weights.}
Quadrature weights do not change through the process, so we fuse them with the transpose  of the gradient operator 
 \begin{align*}
 \DoperTrans_{\MB{\qweights}}
        =
  \Dmat^{\Transpose}
 \MB{\qweights},
 \end{align*}
where $\DoperTrans_{\MB{\qweights}}:\D{R}^{\mandeld\NQ} \rightarrow \D{R}^{d\NI}$ can be interpreted as a weighted discrete divergence operator.

\begin{algorithm}
\setstretch{1.5}
\caption{Pseudo-algorithm of the displacement-based Newton-PCG solver}\label{alg:main}
 \begin{algorithmic}[1]

     \State{\textbf{Initialize: }}
    \State   $\MB{\Dperdisp}_{(0)}$, $\MB{e}$
                \Comment{initial displacement, macroscopic strain}
     \State  $ \eta^{\text{NW}}$, $\eta^{\text{CG}}$ \Comment{Newton- and CG-tolerance}
      \State{$it_\text{max}^{\text{NW}}$, $it_\text{max}^{\text{CG}}$} \Comment{max. iterations Newton and CG}
  \State
         \For{$i = 0, 1, 2, \dots, it_\text{max}^{\text{NW}} $} \Comment{Newton iteration}
\State $\MB{\vek{g}}_{(i)}=\dots$\Comment{update internal parameters}
    \State   $\MB{\RHS}_{(i)}=-
        \DoperTrans_{\MB{\qweights}}
        \stress(\MB{e}+
        \Doper \MB{\Dperdisp}_{(i)}, \MB{\vek{g}}_{(i)})
        $ \Comment{right-hand side}

            \State  $\Dmaterial_{(i)}
                     =
                    \dfrac{\partial \stress}{\partial \strain}(\MB{e}+
                     \Doper \MB{\Dperdisp}_{(i)})
                    $
                    \Comment{material tangent}       
                   
            \State  Assembly $  \KrefinvF
                    $ 
                    \Comment{Preconditioner assembly - Algorithm \ref{alg:precon}}

   \State{Solve for $\delta\Dperdisp_{(i+1)}$ with \textbf{PCG}:}
    \State{$\qquad \linsysMat_{(i)}\delta\Dperdisp_{(i+1)}
         =
       \MB{\RHS}_{(i)}$  \text{with preconditioner}   $\Krefinv$ \text{in} $it_\text{max}^{\text{CG}}$ \text{steps to accuracy} $\eta^\text{CG}$.} 
            
            \State  $
                    \MB{\Dperdisp}_{(i+1)}
                     =
                    \MB{\Dperdisp}_{(i)}+\delta\MB{\Dperdisp}_{(i+1)}
                    $
                    \Comment{iterative update}
             
           \If{$ \norm{\delta\MB{\Dperdisp}_{(i+1)}} \leq   \eta^{\text{NW}} \norm{\MB{\Dperdisp}_{(i+1)}}$} 
    \State{Proceed to line~\ref{line:return}} \Comment{\text{Newton's method converged}}
    \EndIf          
         \EndFor{}
          \State {\textbf{return} $\MB{\Dperdisp}_{(i+1)}$\label{line:return}}
\end{algorithmic}
\end{algorithm}
\subsection{Assembly of $\KrefF$}
It may be useful to reassemble the preconditioner with updated~$\rmaterial_{(i)}$, whenever~$\material_{(i)}$ significantly changes with respect to the previous Newton step, with $\material_{(i-1)}$. 
However, the use of matrix-free operators $\Doper, \DoperTrans_{\MB{\qweights}}, \FourierTrans$ and $\iFourierTrans$ prohibits the direct assembly of $\KrefF$ through matrices, like in~\eqref{eq:M}. Thus, we suggest an efficient algorithm for the assembly of $\KrefF$, that is outlined in Algorithm~\ref{alg:precon}.

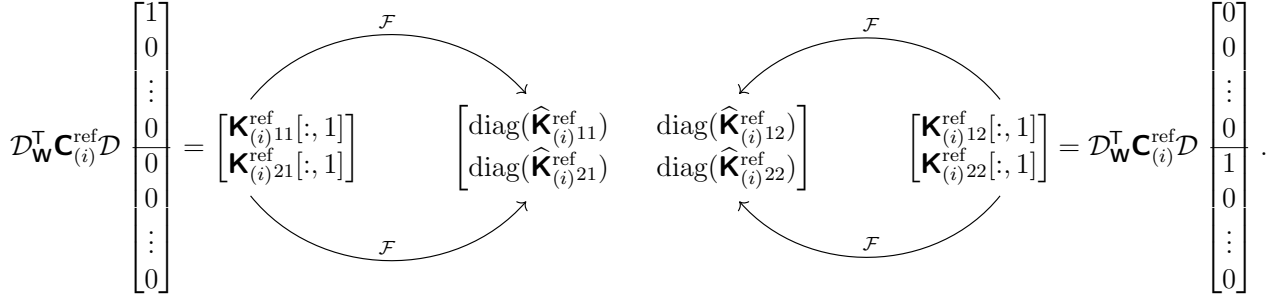
\begin{figure}
\centering
\begin{center}
\adjustbox{scale=0.9,center}{
\begin{tikzcd}
   \DoperTrans_{\MB{\qweights}} \Drmaterial_{(i)} \Doper
 \begin{bmatrix}
   1\\
 0\\
    \vdots\\
    0\\
    \hline
      0\\
       0\\
    \vdots\\
    0\\
   \end{bmatrix}
   \hspace{-7ex}
   &
   =
    \begin{bmatrix}
    \Krefcomp{1}{1}[:,1]  \\
      \Krefcomp{2}{1}[:,1] 
   \end{bmatrix}
    \arrow[r, bend left=50, xshift=-5ex, "\FourierTrans"]
    \arrow[r, bend right=50, xshift=-5ex , "\FourierTrans"]
    &
   \begin{bmatrix}
  \text{diag}( \KrefFcomp{1}{1})& \text{diag}(\KrefFcomp{1}{2})    \\
   \text{diag}(\KrefFcomp{2}{1}) &   \text{diag}(\KrefFcomp{2}{2})
   \end{bmatrix}
 &
   \arrow[l, bend left=50, xshift=5ex,swap,"\FourierTrans"]
    \arrow[l, bend right=50, xshift=5ex ,swap, "\FourierTrans"]  
    \begin{bmatrix}
    \Krefcomp{1}{2}[:,1]  \\
      \Krefcomp{2}{2}[:,1] 
   \end{bmatrix}
   &
   \hspace{-7ex}
   =
       \DoperTrans_{\MB{\qweights}} \Drmaterial_{(i)} \Doper
 \begin{bmatrix}
  0\\
 0\\
    \vdots\\
    0\\
    \hline
      1\\
       0\\
    \vdots\\
    0\\
   \end{bmatrix}.
   \end{tikzcd}
}  
\end{center}
\caption{The schematic procedure of matrix-free assembly of $\KrefF$ for $d\Nn=2$. First columns of blocks $\Krefcomp{1}{1}$ and $\Krefcomp{2}{1}$ are obtained as a result of the matrix-free action of $\Kref$ on the unit impulse vector $\MB{i}^1$. Diagonals $\text{diag}(\KrefFcomp{1}{1})$ and $\text{diag}(\KrefFcomp{2}{1})$ are then computed through $d$-dimensional FFT of $\Krefcomp{1}{1}[:,1]$ and $\Krefcomp{2}{1}[:,1]$, respectively. By analogy, columns of blocks $\Krefcomp{1}{2}$ and $\Krefcomp{2}{2}$ are obtained by the matrix-free action of $\Kref$ on the unit impulse vector $\MB{i}^p$ where $p= (2-1)\Np+1 $. }
\label{fig:diagram}
\end{figure}

First, take a look at (block-periodic) $\bar{\alp}\bar{\beta}$-block $\Krefab\in \D{R}^{\Np\times \Np}$ of $\Kref\in \D{R}^{d\Nn\Np\times d\Nn\Np}$. Thanks to the convolution theorem, the whole diagonal $\text{diag}(\KrefFab)\in \D{R}^{\Np}$ can be obtained by the FFT of any, say the first, row or, because of the symmetry, column of $\Krefab$,
\begin{align*}
    \text{diag}(\KrefFab)={\FourierTrans (\Krefab[1,:])}^{\Transpose}=\FourierTrans (\Krefab[:,1])
\end{align*}
where a colon indicates a complete column or row.
Before the FFTs, we have to compute one column $\Krefab[1,:]$ for each of $(d\Nn)^2$ blocks $\Krefab$ of $\Kref$.  
Consider a unit impulse vector $\MB{i}^p \in \D{R}^{d\Nn\Np}$ that has only one non-zero element equal to $1$ on the $p$-th position. When we apply $\Kref$ to vector $\MB{i}^1$, we obtain the first columns of $d\Nn$ blocks $\Krefcomp{\bar{\alp}}{1}$. From the structure of $\Kref$ visible in \eqref{eq:M} it is obvious that we need $d\Nn$ vectors $\MB{i}^p$ to obtain all $(d\Nn)^2$ columns $\Krefab[1,:]$, where $p= (\bar{\beta}-1)\Np+1 $ and $\bar{\beta} \in \{ 1, \dots, d \Nn \} $. The whole procedure is schematically shown in Fig.~\ref{fig:diagram}.
%

\subsection{Pseudo-Inverse of $\KrefF$}
Once we have all diagonal blocks $\KrefFab$,
we may proceed to the computation of the pseudo-inverse of $\KrefF$. By a proper row and column reordering, it can be seen that the pseudo-inverse of the block diagonal matrix $\KrefF$ is equivalent to the pseudo-inverse of $\Np$ (number of pixels/stencils) submatrices 
   \begin{align}\label{eq:block_inversion}
  \begin{bmatrix}
   \KrefFcomp{1}{1} [J,J]& \dots & \KrefFcomp{1}{\bar{\beta}}   [J,J] \\
   \vdots& \ddots & \vdots  \\
      \KrefFcomp{\bar{\alp}}{1}  [J,J] &  \dots & \KrefFcomp{\bar{\alp}}{\bar{\beta}}[J,J]  
   \end{bmatrix}^{-1} 
   \in \D{R}^{d\Nn\times d\Nn}, \quad \text{where} \quad J\in \{1, \dots, \Np \}.
 \end{align}
The $(\Np-1)$ submatrices are of full rank and thus directly invertible. Only one submatrix, corresponding to the zero frequency Fourier mode, is singular and has to be treated separately. This block has exactly $d$ null eigenvalues corresponding to $d$ rigid-body modes. We compute the (Moore-Penrose) pseudo-inverse of this block instead of its inversion\footnote{Please note that the Moore-Penrose pseudo-inverse is depicted by $\dagger$ in Algorithm~\ref{alg:precon}.}. The pseudo-inverse can be computed exactly by restriction onto the space orthogonal to the kernel of the singular block. For any specific type of FE and the corresponding discretisation stencil, the kernel can be exactly identified.

\begin{algorithm}
\setstretch{1.5}
\caption{Pseudo-algorithm of reference material based preconditioner assembly}\label{alg:precon}
 \begin{algorithmic}[1]
 
     \State{\textbf{Initialize: }}
    \State  $\Drmaterial_{(i)}$
                \Comment{spatially uniform reference material}

        \State           
    \For{$ \bar{\beta} = 1,\dots, d\Nn$} \Comment{loop over $d$ vectors} 
        \State $p= (\bar{\beta}-1)\Np+1$
                 \Comment{column index}
        \State $\vek{c}_{\bar{\beta}}=\DoperTrans_{\MB{\qweights}} \Drmaterial_{(i)} \Doper \MB{i}^{p}$
                \Comment{$p $-th column of $\Kref$}
        \For{$\bar{\alp} = 1,\dots, d\Nn$}
            \State   $\text{diag}(\KrefFab)=\FourierTrans( \vek{c}_{\bar{\beta}}[(\bar{\alp}-1)\Np+1:\bar{\alp}\Np]) $                    \Comment{assign to $\KrefFab$  diagonals}
    
        \EndFor
    \EndFor
        \State
            \Comment{pseudo-inverse of singular submatrix of $\KrefF$}
         
        \State
            $\KrefinvabF[1,1]
            =
             \begin{bmatrix}
   \KrefFcomp{1}{1} [1,1]& \dots & \KrefFcomp{1}{\bar{\delta}}   [1,1] \\
   \vdots& \ddots & \vdots  \\
      \KrefFcomp{\bar{\gamma}}{1}  [1,1] &  \dots & \KrefFcomp{\bar{\gamma}}{\bar{\delta}}[1,1]  
   \end{bmatrix}^{\dagger}_{\bar{\alp}\bar{\beta}}$
            \Comment{$1$-th block}
    
    \For{$J = 2,\dots, \Np $}
            \Comment{inverse of remaining submatrices of $\KrefF$}
         
        \State
            $\KrefinvabF[J,J]
            =
             \begin{bmatrix}
   \KrefFcomp{1}{1} [J,J]& \dots & \KrefFcomp{1}{\bar{\delta}}   [J,J] \\
   \vdots& \ddots & \vdots  \\
      \KrefFcomp{\bar{\gamma}}{1}  [J,J] &  \dots & \KrefFcomp{\bar{\gamma}}{\bar{\delta}}[J,J]  
   \end{bmatrix}^{-1}_{\bar{\alp}\bar{\beta}}$
            \Comment{$J$-th block}
    \EndFor
\end{algorithmic}
\end{algorithm}

 \section{Numerical Experiments}\label{sec:experiments}
 We demonstrate the numerical behavior of the proposed approach on several examples. In general, we compare our displacement-based (DB) FE scheme, described in the previous sections, with the CG accelerated strain-based (SB) Fourier-Galerkin method taken from \cite{VoZeMa2014FFTH,ZeGeVoPeGe2017,DEGEUS2017412}. All results were obtained with the $\mu$Spectre software, an open-source platform for efficient FFT-based continuum mesoscale modelling, which is freely available at 
 \url{https://gitlab.com/muspectre/muspectre}. The software package includes the examples, which are described in the following sections.

 \emph{Termination criteria.}
To obtain comparable results, we have to choose the corresponding termination criteria for both SB and DB schemes.
The Newton's method stops when the relative norm of strain increment drops below the tolerance $ \eta^{\text{NW}}$,  ${\norm{\delta\MB{\symgradM \Dperdisp}_{(i+1)}}}\leq  \eta^{\text{NW}} {\norm{\MB{\symgradM \Dperdisp}_{(i+1)}}}$.
The CG solver for SB scheme is stopped when the relative norm of the residual drops below the tolerance $ \eta^{\text{CG}}$, $ {\norm{\MB{r}^{\text{SB}}_{(i+1)}}}\leq \eta^{\text{CG}}{\norm{\MB{r}^{\text{SB}}_{(0)}}}$, while the PCG solver for DB schemes is stopped when the relative $\Krefinv$-norm of the residual drops below the tolerance $ \eta^{\text{CG}}$,
 ${\norm{\MB{r}^{\text{DB}}_{(i+1)}}_{\Krefinv }}\leq \eta^{\text{CG}}{\norm{\MB{r}^{\text{DB}}_{(0)}}_{\Krefinv }} $.
 This nontrivial equivalence will be discussed in  detail in an upcoming publication~\cite{Ladecky2022}.

  
\subsection{Linear Steady-State Thermal Conduction Problem}\label{sec:example_scalar}
In the first example, we demonstrate the oscillation-free character of gradient fields arising from the FE discretisation. For this purpose, we reconstruct the benchmark problem from \cite[Section~3.7.1]{Leuschner2018} or \cite[Section~3.2]{Brisard2010FFT}, where the Fourier-Galerkin methods exhibit significant discretisation artifacts.

\begin{figure}
  \centering
  \includegraphics[width=0.255\textwidth]{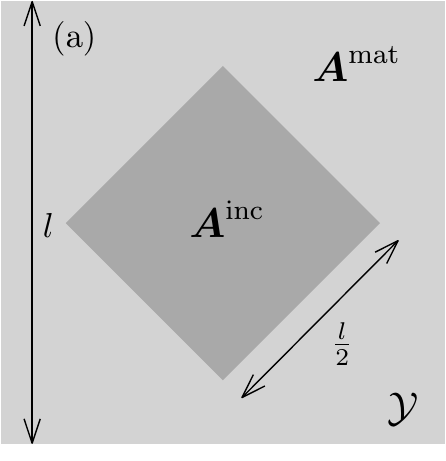}
  \includegraphics[width=0.35\textwidth]{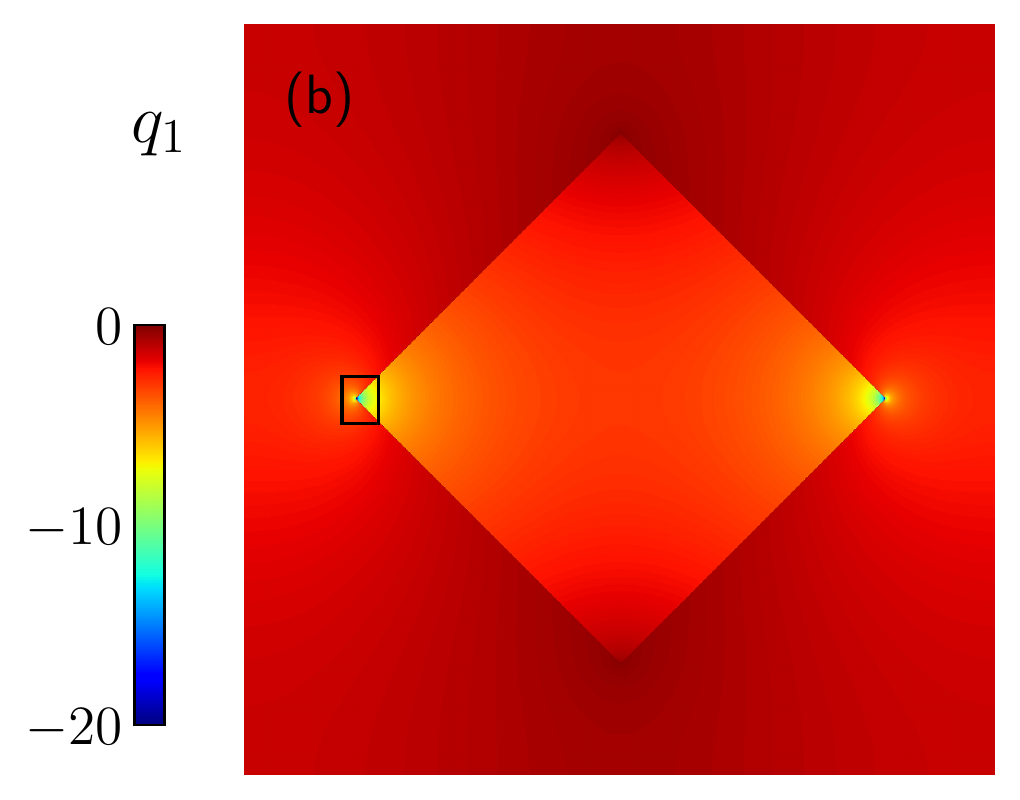}
  \includegraphics[width=0.34\textwidth]{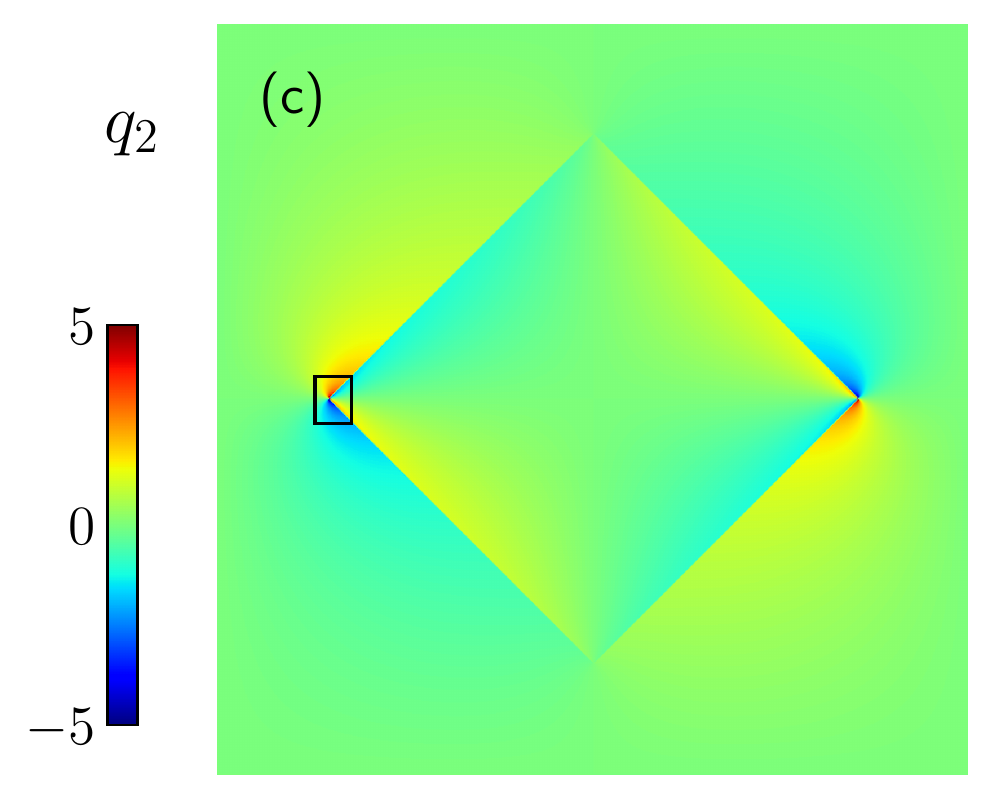}
  \caption{A linear heat transfer problem from Section \ref{sec:example_scalar}. The square periodic unit cell $\puc$ with a square inclusion \textbf{(a)}. The flux field component $q_{1}$ \textbf{(b)} and $q_{2}$ \textbf{(c)} arising from average temperature gradient $\macrostrain = \left[
   0.01  ,0.0 \right]^{\Transpose}$. Results are obtained with one-node FE stencil ($\Nn= 1$) with two linear triangular elements discretisation and $815$ nodes in both directions ($\NI= 815^2$).}
  \label{fig:diamond_puc}
\end{figure}

\begin{figure}[htbp]
  \centering
 \includegraphics[width=1\textwidth]{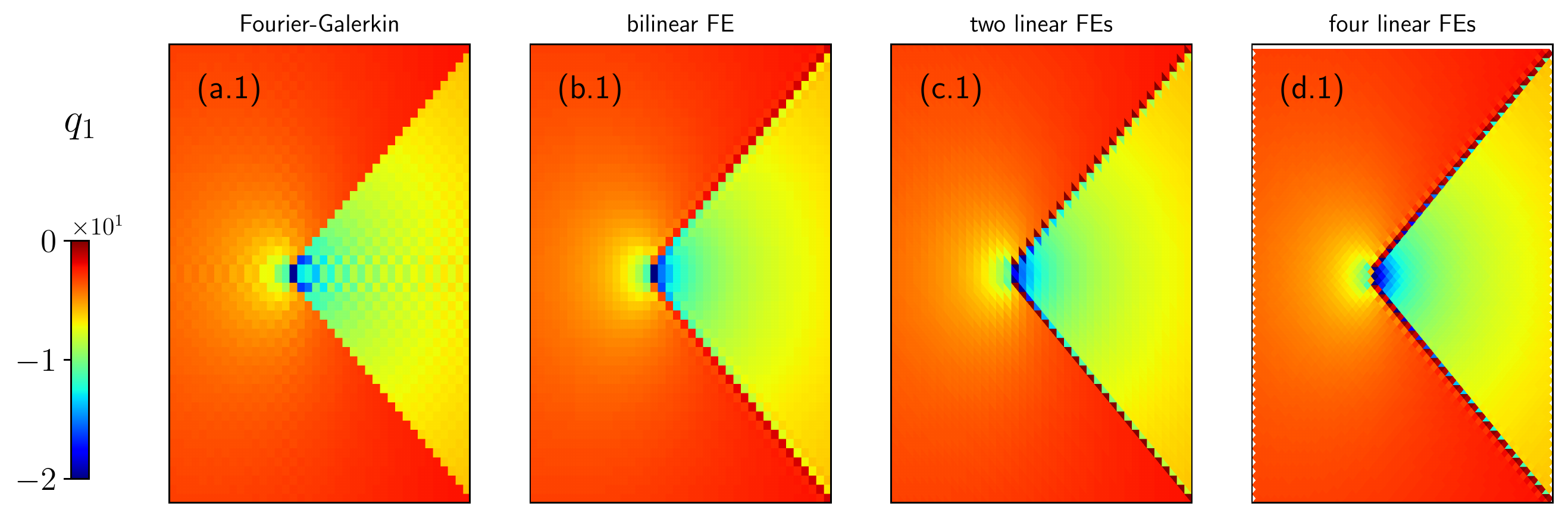}
 \includegraphics[width=1\textwidth]{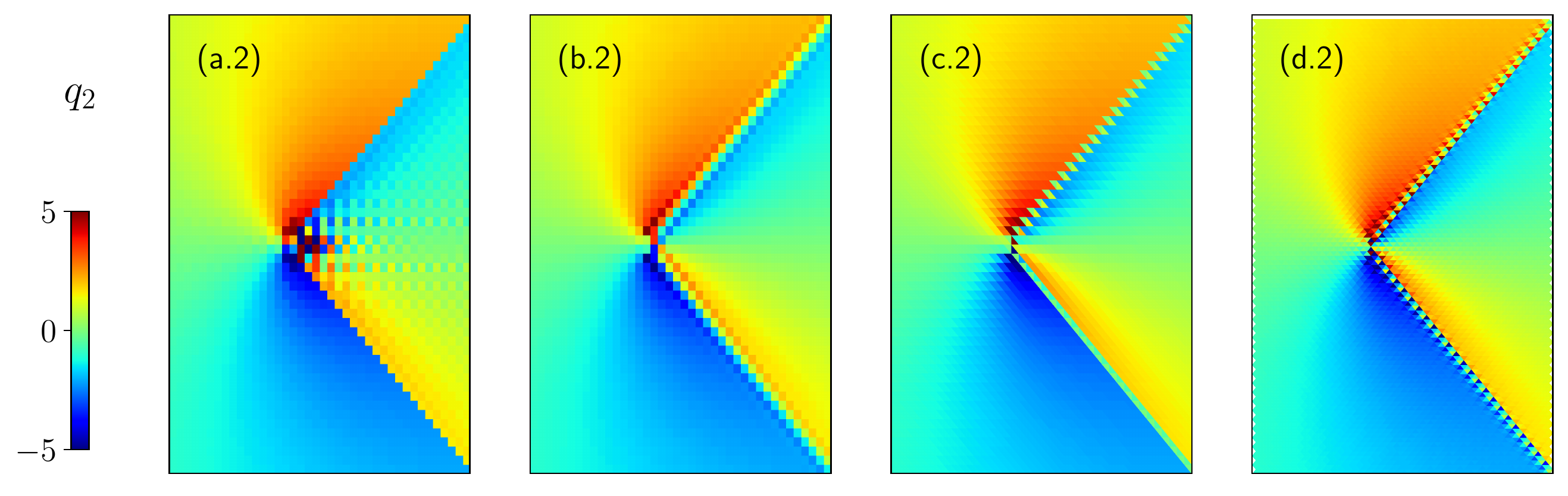}
  \caption{Local heat flux field components $q_1$ \textbf{(1)} and $q_2$ \textbf{(2)} from experiment in Section \ref{sec:example_scalar}, obtained with the Fourier-Galerkin method is shown in the column \textbf{(a)}, one-node FE stencil ($\Nn= 1$) with one bilinear rectangular element is shown in the column \textbf{(b)}, one-node FE stencil ($\Nn= 1$) with two linear triangular elements is shown in the column \textbf{(c)}, and two-node FE stencil ($\Nn= 2$) with four linear triangular elements is shown in the column \textbf{(d)}.}
  \label{fig:diamond_q1}
\end{figure}

\begin{figure}[htbp]
  \centering
 \includegraphics[width=1\textwidth]{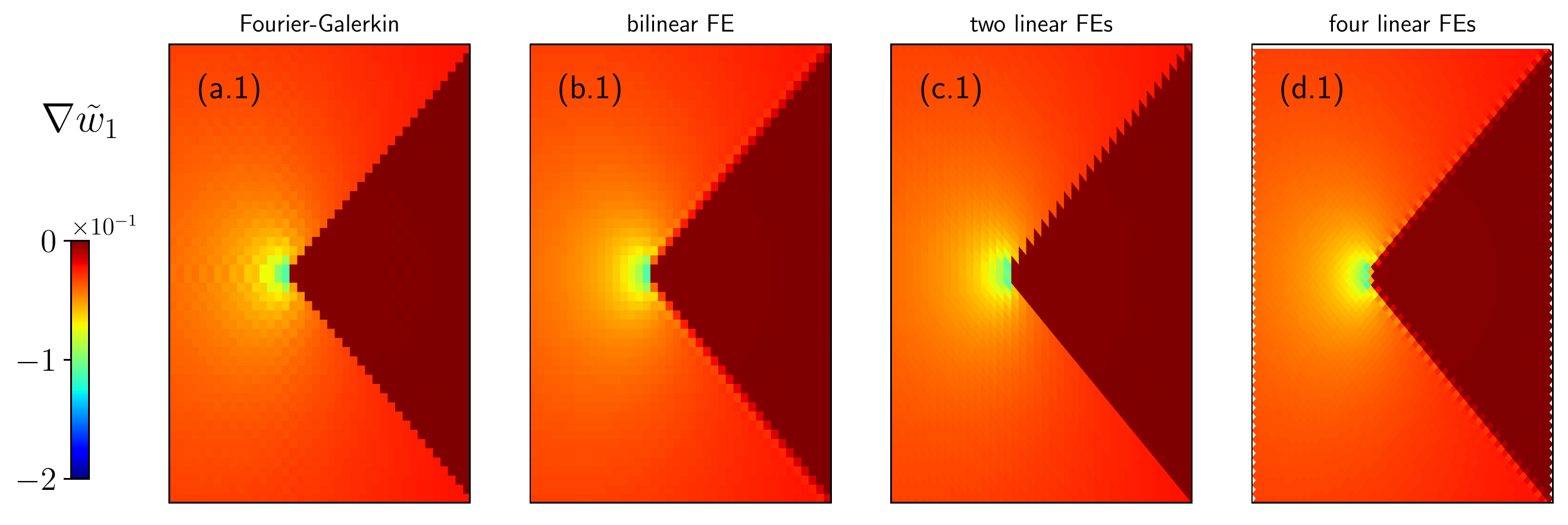}
 \includegraphics[width=1\textwidth]{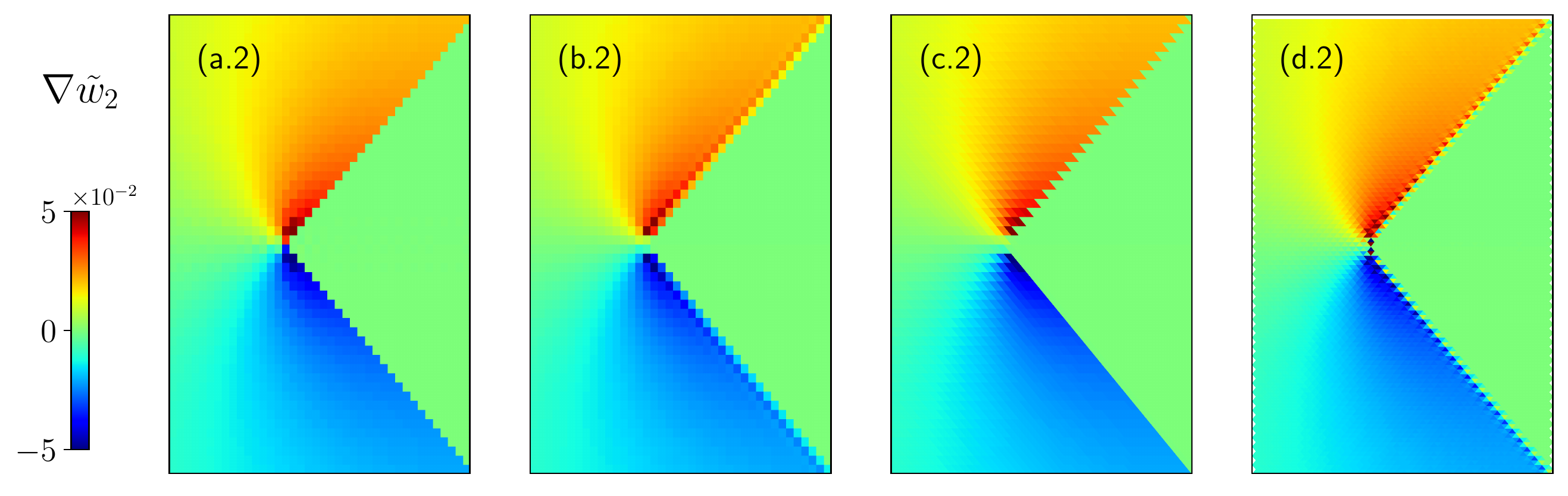}
  \caption{Local temperature gradient field components $\grad \pertemp_1$ \textbf{(1)} and $\grad \pertemp_2$ \textbf{(2)} from experiment \ref{sec:example_scalar}, obtained with the Fourier-Galerkin method is shown in the  column \textbf{(a)}, one-node FE stencil ($\Nn= 1$) with one bilinear rectangular element is shown in the  column \textbf{(b)}, one-node FE stencil ($\Nn= 1$) with two linear triangular elements is shown in the column \textbf{(c)}, and two-node FE stencil ($\Nn= 2$) with four linear triangular elements is shown in the  column \textbf{(d)}.}
  \label{fig:diamond_g1}
\end{figure}

We consider a scalar problem of linear heat transfer, where we look for the flux field $\flux$ satisfying the weak balance condition \eqref{eq:weak_form_scalar};  see~\ref{sec:appendix_termal} for more details. The microstructure is defined by the square periodic unit cell $\puc$, as sketched on the left-hand side of Fig.~\ref{fig:diamond_puc}. The composite microstructure consists of an insulating matrix with the conductivity $\MB{\conductmaterial}^{\text{mat}}=100\, \mat{\identity}$, and a conducting inclusion with the conductivity $\MB{\conductmaterial}^{\text{inc}}=100\, \MB{\conductmaterial}^{\text{mat}}$. 
An average temperature gradient $\macrostrain = \left[
   0.01  ,0.0 \right]^{\Transpose}$ is applied. The number of pixels is $815^{2}$, and the material coefficients are constant per pixel.
Components of the global flux field $\flux$ are shown in Fig.~\ref{fig:diamond_puc};  $q_{1}$ in the middle and $q_{2}$ on the right-hand side. The regions of details  depicted in Fig.~\ref{fig:diamond_q1} and Fig.~\ref{fig:diamond_g1} are highlighted by the black rectangles in Fig.~\ref{fig:diamond_puc}. 

In Fig.~\ref{fig:diamond_q1}, we show the details of heat fluxes for various discretisations: the Fourier-Galerkin method in the column (a), the one-node FE stencil ($\Nn= 1$) with one bilinear rectangular element (Fig.~\ref{fig:tiling} (a)) in the column (b), the one-node FE stencil ($\Nn= 1$) with two linear triangular elements (Fig.~\ref{fig:tiling} (b)) in the column (c), and the two-node FE stencil ($\Nn= 2$) with four linear triangular elements (Fig.~\ref{fig:tiling} (c)) in the column (d).

The Fourier-Galerkin method exhibits strong oscillations through the region, while all FE solutions are devoid of oscillations in the interior of the domains occupied by a single phase as discretisation discrepancies remain confined to the vicinity of the phase boundaries. For instance, triangular discretisations reduce the phase boundary discretisation artifacts to the two pixel-wide layer around the phase boundary.

The zigzag patterns on the phase boundary arise from the pixel-based geometry. 
If the elements can capture the interface of the two phases exactly, we do not get any discretisation artifacts, as can be seen in Fig.~\ref{fig:diamond_q1} and Fig.~\ref{fig:diamond_g1} in the  column (c). This speaks in favour of using FE over the Fourier-Galerkin discretisation.

\subsection{Small-Strain Elasticity Problem}\label{sec:example_small_elasticity}
The second example focuses on the effect of the preconditioner on the number of PCG iterations with respect to the number of discretisation nodes $\NI$ and phase contrast $\rho$.  For this purpose, we use Hashin's coated sphere construction adapted from \cite[Section~4.1]{Schneider2016} and the references therein.

\begin{figure}
  \centering
  \includegraphics[width=1.\textwidth]{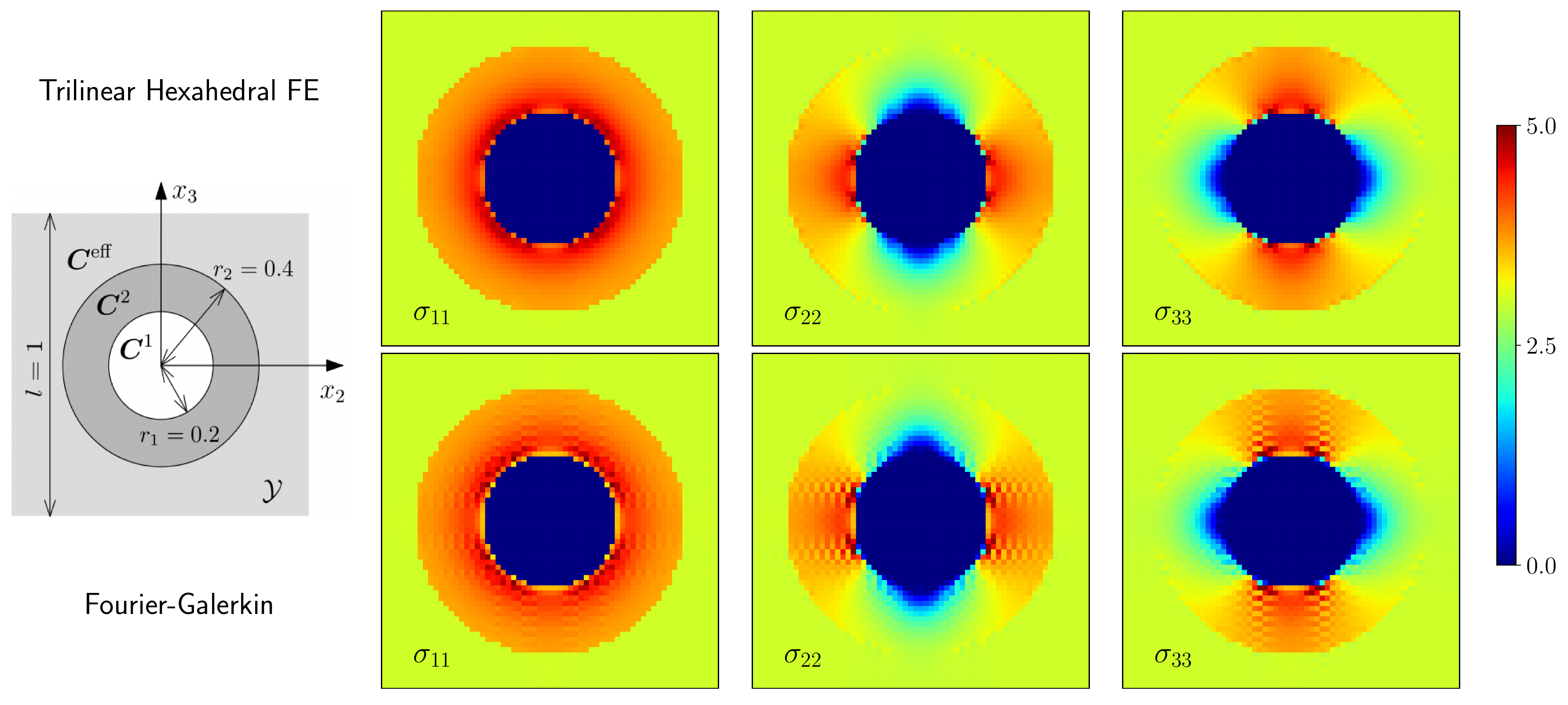}
  \caption{Two-dimensional sections at $x_1 =0.5 $ of the cubic periodic unit cell $\puc$ with a coated sphere inclusion. Radii $r_{1} = 0.2$, $r_{2} = 0.4$ and the domain size $l=1$ (left). Components of the local stress fields $\sigma_{\alpha\alpha}$ for trilinear FE discretisation (top) and Fourier-Galerkin discretisation (bottom) with the number discretisation nodes $\NI=65^3$ .}
  \label{fig:sphere_puc}
\end{figure}

We choose a linear small-strain elastic problem described in Section \ref{sec:small_strain_elasticity}. The three-phase microstructure representing a coated sphere in the matrix with effective material properties is depicted in Fig.~\ref{fig:sphere_puc}, with the core radius $r_{1} = 0.2$, annulus-shaped coating outer radius $r_{2} = 0.4$ and the cubic domain edge length $l=1$. An~average macroscopic strain $\macrostrain =\left[
   1,0  ,0,0 ,0 ,0  \right]^{\Transpose}$ is applied. 
We assume isotropic phases with bulk and shear moduli $\text{K}_1,\text{G}_1$ in the core, $\text{K}_2,\text{G}_2$ in the coating and $\text{K}_\text{eff},\text{G}_\text{eff}$ in the surrounding matrix. The bulk moduli $\text{K}_1, \text{K}_2$ are chosen in a way that resulting response of the unit cell is equivalent to the response of a homogeneous material with $\text{K}_\text{eff}$. As a consequence, the bulk moduli $\text{K}_1, \text{K}_2$ have to be balanced for particular phase contrast $\rho =  {\text{K}_2}/{\text{K}_1}$.

First, in accordance with \citet[Section~4.1.3]{Schneider2016}, we set $\rho = 10^{3} $ and the remaining parameters to
\begin{align*}
    \text{K}_1 \doteq 0.00132060, \quad \text{K}_2 \doteq 1.3206033, \quad \text{K}_\text{eff} \doteq 1.0,\\
   \text{G}_1 \doteq  0.00079236, \quad \text{G}_2 \doteq  0.7923620, \quad \text{G}_\text{eff} \doteq 0.6.
    \end{align*}
 Two-dimensional sections at $x_1 =0.5 $ of global stress field components are shown in Fig.~\ref{fig:sphere_puc} right. Trilinear FE discretisation (the top row in Fig.~\ref{fig:sphere_puc}) generates oscillation free results compared to the Fourier-Galerkin discretisation (the bottom row in Fig.~\ref{fig:sphere_puc}).

Second, we are interested in how our preconditioned scheme behaves with respect to the number of discretisation nodes $\NI$ and varying phase contrast $\rho $. The convergence of PCG depends on the choice of reference material $\rmaterial$. We compare two cases: the first $\rmaterial_{\mat{\identity_{s}}}= \mat{\identity_{s}}$, with $ \mat{\identity_{s}} \in  \D{R}^{\mandeld \times \mandeld} $ being the symmetrized identity tensor $({I_{\mathrm{s}}})_{\alpha\beta\iota\kappa} =\frac{1}{2}(\delta_{\alpha\iota}\delta_{\beta\kappa}+\delta_{\alpha\kappa}\delta_{\beta\iota})$ in the Mandel notation, and secondly $\rmaterial_{\text{mean}}= \frac{1}{| \puc |}\sum_{Q=1}^{\NQ} \material(\vek{x}_{\mathrm{q}}^Q) w^Q$, where  $\rmaterial_{\text{mean}}$ is the mean stiffness matrices over $\puc$.

\begin{figure}[htbp]
  \centering
 \includegraphics[width=0.95\textwidth]{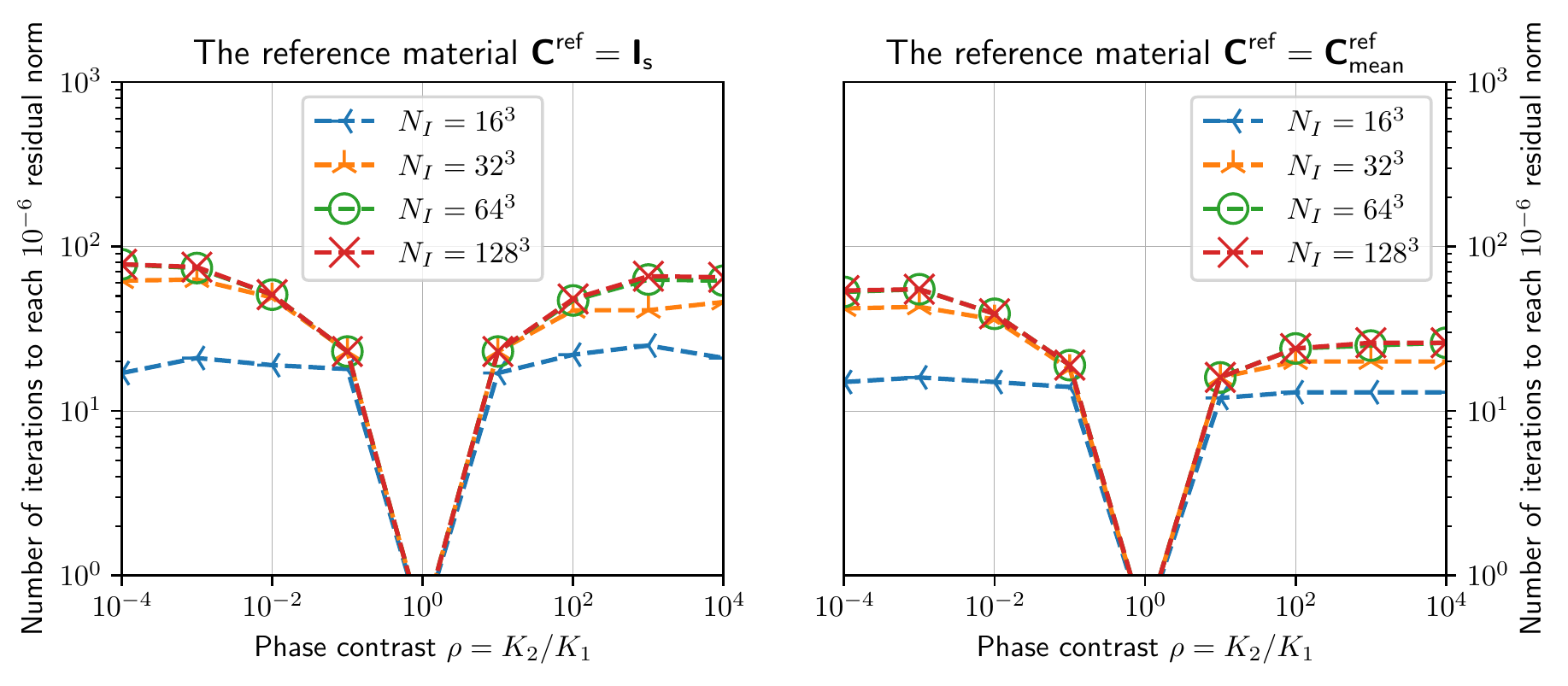}
  \caption{The number of PCG iterations for trilinear FE discretisation for different phase contrasts $\rho$ and number of discretisation nodes $\NI$. Termination parameter for linear solver $\eta^{\text{CG}}=10^{-6}$.}
  \label{fig:iterations}
\end{figure}

The preconditioner with mean reference material $\rmaterial_{\text{mean}}$ exhibits better performance in all studied cases, see Fig.~\ref{fig:iterations}. The numbers of iterations slowly increases with the growing $\NI$ until it stabilizes for sufficiently fine discretisations. In addition, $\rmaterial_{\text{mean}}$ significantly reduced the phase contrast sensibility, especially for $\rho>1$ (softer sphere core).

\subsection{Finite Strain Elasto-Plastic Problem}\label{sec:elastoplastic}
The purpose of the last example is twofold. First, we demonstrate the applicability of the approach to real-world problems in the finite strain setting, and the effect of nonphysical oscillations on the results. Second, we point out the equivalence of our DB FE scheme with SB scheme with FE projection operator recently proposed by \citet{LeuteR2021}. The equivalence of these two approaches is briefly explained later in Section~\ref{sec:comparison_gradients}. 


To this purpose, we adapt the example from Section 5.5 of~\citet{DEGEUS2017412}. The example studies a sample of a dual-phase steel obtained by a scanning electron microscope. Responses of the material phases are elastic and homogeneous in the elastic part of deformation with Young's moduli $E=200\, \text{GPa}$ and Poisson's ratios $\nu=0.3$, and differ in the parameters of linear hardening in the plastic regions, see~\cite[Section 5]{DEGEUS2017412} for more details on the material model. 

The yield stress $\tau_{y}=\tau_{y0}+H\varepsilon_{\text{p}}$ evolves with respect to plastic strain $\varepsilon_{\text{p}}$, initial yield stresses $    \tau^{\text{hard}}_{y0}=2\tau^{\text{soft}}_{y0}=0.003E $, and hardening moduli
$    H^{\text{hard}}_{0}= 2H^{\text{soft}}_{0}=0.01 E$. Total macroscopic deformation gradient
\begin{align}\label{eq:F}
\mat{F}=\frac{\sqrt{3}}{2} \footnotesize{\begin{bmatrix}
   0.995 & 0 \\
   0 & -0.995
\end{bmatrix}}
\end{align} is applied in $5$~load increments. 

We solved this problem with the following schemes: the  standard SB Fourier-Galerkin scheme with Fourier projection operator from \cite{ZeGeVoPeGe2017,DEGEUS2017412}, SB scheme with two linear triangular FEs and the FE projection operator from \cite{LeuteR2021}, the DB FE scheme with two linear triangular FEs, and the DB FE scheme with one bilinear rectangular FEs. The Newton tolerance $\eta^{\text{NW}}=10^{-5}$ and (P)CG tolerance $\eta^{\text{CG}}=10^{-5}$. 

\begin{figure}[htbp]
  \centering
 \includegraphics[width=1\textwidth]{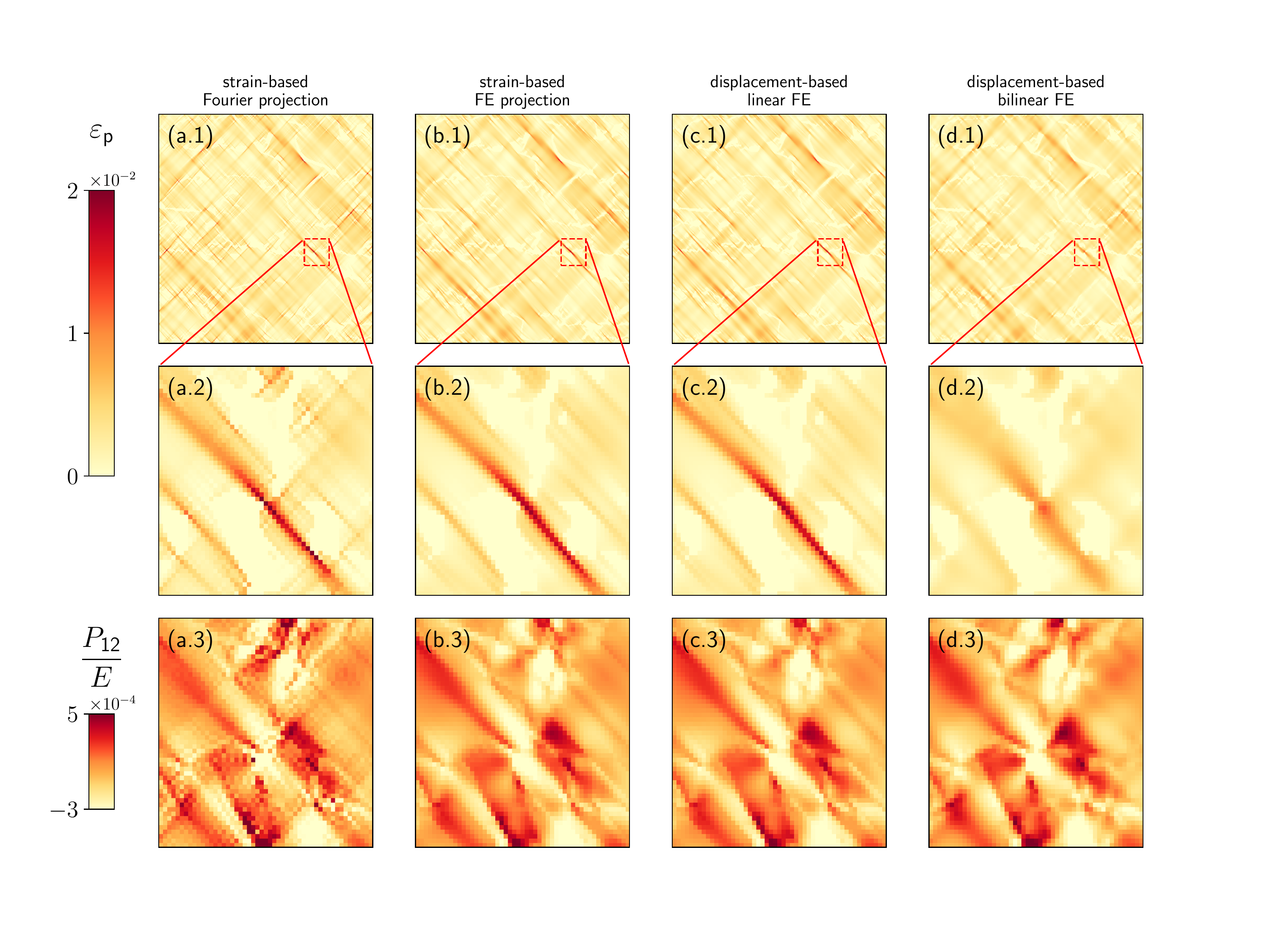}
  \caption{Global plastic strains $\varepsilon_{\text{p}}$ in dual-phase steel with applied deformation gradient~\eqref{eq:F} in row \textbf{(1)} with local detials in row \textbf{(2)}. Row \textbf{(3)} shows accompanying normalized shear stresses $P_{12}$ in detailed area.
  Discretisation schemes in columns: \textbf{(a)} the  standard SB scheme with Fourier projection operator, \textbf{(b)} the SB scheme with FE projection operator with two linear triangular elements, \textbf{(c)}  the DB FE scheme with two linear triangular elements, and \textbf{(d)} the DB FE scheme with one bilinear rectangular elements. All quantities are averaged per pixel.  
  }
  \label{fig:plastic_global}
\end{figure}

First, the distributions of global plastic strain $\varepsilon_{\text{p}}$ obtained for these four approaches are shown in the first row of Fig.~\ref{fig:plastic_global}.
The regions of details (the second row) uncover the checkerboard patterns in the plastic strain fields of the  standard SB Fourier-Galerkin solution (a.2), that are a direct consequence of the oscillating stress field (a.3). The other three schemes, columns (b) to (d), produce solutions without oscillations.

Second, the number of Newton's method steps and the total number of (P)CG iterations needed to solve the problem with these four approaches are shown in Table~\ref{tab:equvivalence}. The table highlights the equivalence of our DB scheme and the SB scheme presented by~\citet{LeuteR2021}, if equivalent discretisations are used. 

\begin{table}[htbp]
\centering
 \begin{tabular}{m{0.15\textwidth}|>{\centering}m{0.19\textwidth} >{\centering}m{0.15\textwidth}>{\centering}m{0.15\textwidth}>{\centering\arraybackslash}m{0.19\textwidth}} 
  & \multicolumn{2}{c}{strain-based (SB)} &  \multicolumn{2}{c}{displacement-based (DB)}\\[10pt]
 & Fourier projection & FE projection &   linear FE & bilinear FE  \\[10pt]
 \hline\\[-5pt]
Newton~steps & 11 & \textbf{9} & \textbf{9} & 10\\[10pt]
(P)CG~steps & 1013 & \textbf{862} & \textbf{862} & 761 
 \end{tabular}
 \caption{The number of Newton's method steps and the total number of CG or PCG steps required to solve the finite strain elasto-plastic problem of Section~\ref{sec:elastoplastic}, with Newton tolerance $\eta^{\text{NW}}=10^{-5}$ and (P)CG tolerance $\eta^{\text{CG}}=10^{-5}$.
 Discretisation approaches from left to right:  the  standard SB Fourier-Galerkin scheme with Fourier projection operator, SB scheme with FE projection operator with two linear triangular elements, the DB FE scheme with two linear triangular elements, and the DB FE scheme with one bilinear rectangular element per pixel. Numbers in boldface highlight the equivalence of our DB FE scheme and SB FE scheme presented by~\citet{LeuteR2021}.}
 \label{tab:equvivalence}
\end{table}

  \section{Comparison with related FFT-based schemes}\label{sec:comparison} 
Several FFT-based computational homogenization schemes exist~\cite{Schneider2021,Lucarini_2021}. 
 An interested reader may therefore find a comparison and placement of our approach in the context of contemporary literature useful. 

Recall that our approach is derived from the weak form of the mechanical equilibrium condition \eqref{eq:weak_form} with an unknown displacement field. The equilibrium \eqref{eq:weak_form} is discretized in the standard Galerkin manner with the FE basis functions. The nonlinear nodal equilibrium \eqref{eq:non_lin_system} is linearized by the Newton's method, and the system of linear equations \eqref{eq:lin_system} is solved by the PCG method. Favourable convergence property of the PCG method is guaranteed by the reference material based preconditioner \eqref{eq:Kref}, which fast application builds on FFT. 

\subsection{The Connection with Strain-Based Approaches}\label{sec:comparison_gradients}

Unlike the DB FE, most spectral methods use strains (gradients) as unknown. SB approaches, like those in \cite{LeuteR2021,ZeGeVoPeGe2017,Schneider2016}, typically use the projection operator to enforce the compatibility of strain fields. 
To reveal a link between the DB and SB approaches, recall the preconditioned scheme \eqref{eq:prec_scheme},
 \begin{align}\label{eq:precondit_chceme}
        \underbrace{
        (\Dmat^{\Transpose} \MB{\qweights} \Drmaterial_{(i)} \Dmat )^{-1}}_{\Krefinv}
         \underbrace{
         \Dmat^{\Transpose}
          \MB{\qweights}
        \MB{C}_{(i)}
         \Dmat
          }_{\linsysMat_{(i)}}
         \delta\Dperdisp_{(i+1)}
         =&
         -
         \underbrace{
        (\Dmat^{\Transpose} \MB{\qweights} \Drmaterial_{(i)} \Dmat )^{-1}}_{\Krefinv}
          \underbrace{
         \Dmat^{\Transpose} \MB{\qweights} \stress(\MB{e}+\Dmat \Dperdisp_{(i)}, \MB{\vek{g}}_{(i)})
         }_{-\RHS_{(i)}},
\end{align}
 where we omit the FFTs for simplicity.
In the case of linear triangles or tetrahedral elements with a single quadrature point per element, all quadrature weights $w^Q$ are equal. Then the multiplication by quadrature weights $\MB{\qweights}$ can be left out in \eqref{eq:precondit_chceme}, leading to
 \begin{align*}
        \underbrace{
        (\Dmat^{\Transpose} \Drmaterial_{(i)} \Dmat )^{-1}}_{\Krefinv}
         \underbrace{
         \Dmat^{\Transpose}
        \MB{C}_{(i)}
         \Dmat
          }_{\linsysMat_{(i)}}
         \delta\Dperdisp_{(i+1)}
         =&
         -
         \underbrace{
        (\Dmat^{\Transpose} \Drmaterial_{(i)} \Dmat )^{-1}}_{\Krefinv}
          \underbrace{
         \Dmat^{\Transpose}\stress(\MB{e}+\Dmat \Dperdisp_{(i)}, \MB{\vek{g}}_{(i)})
         }_{-\RHS_{(i)}}.
\end{align*}
Next, we replace the iterated unknown $\Dperdisp_{(i)}$ with its gradient $\symgradM \Dperdisp_{(i)}$, recognizing that $\symgradM \Dperdisp_{(i)} = \Dmat\MB{\Dperdisp}_{(i)}$. 
After the multiplication with $\Dmat$ from the left-hand side, we finally obtain
    \begin{align}\label{eq:projection_chceme}
     \underbrace{\Dmat (\Dmat^{\Transpose} \Drmaterial_{(i)} \Dmat )^{-1}
         \Dmat^{\Transpose}}_{\MB{\Gamma}_{(i)}^{0} }
        \MB{C}_{(i)}
          \delta \symgradM \perdisp_{(i+1)}
         =&
         -
         \underbrace{
         \Dmat(\Dmat^{\Transpose}\Drmaterial_{(i)}\Dmat )^{-1}
         \Dmat^{\Transpose}
         }_{\MB{\Gamma}_{(i)}^{0} }
         \stress(\MB{e}+\symgradM \perdisp_{(i)}, \MB{\vek{g}}_{(i)})
        ,
\end{align}
where $\MB{\Gamma}_{(i)}^{0}:\D{R}^{\mandeld\NQ}\rightarrow \D{R}^{\mandeld\NQ}$ stands for the discretized periodic Green's operator. Leute et al.~\cite{LeuteR2021} described that by setting $\Drmaterial_{(i)} =\MB{I_s} $, $\MB{\Gamma}_{(i)}^{0}$ projects  an  arbitrary  field from $\D{R}^{\mandeld\NQ}$ to its closest compatible  part in the least square sense  with respect to the $L^2$-norm.

Therefore, this section demonstrates that the schemes \eqref{eq:precondit_chceme} and \eqref{eq:projection_chceme} are equivalent and generate equivalent solutions in every step of the CG in exact arithmetic. If corresponding stopping criteria are used, CG yields the same approximate solutions.
Thus, the only decision-making argument is the possibility of efficient implementation. 
  
\subsection{The Connection with FEM-FFT Approaches}
 To the best of our knowledge, our method shares the most similarities with the linear hexahedral elements (FFT-$Q_{1}$ Hex) formulation by Schneider at al.~\cite{schneider_fft-based_2017} and Fourier-Accelerated Nodal Solver (FANS) by Leuschner and Fritze \cite{Leuschner2018}. The novelty of our approach lies in the following:
\begin{itemize}
    \item \emph{The gradient operator.} 
    Similarly to FFT-$Q_{1}$ Hex  and FANS, the gradient field is derived with respect to the FE approximation. However, we do not express the discrete gradient operator $\Dmat$ in the Fourier space, but keep it in the real space. The direct convolution with a short gradient kernel is cheaper than the Fourier convolution via forward and inverse FFTs. We use the Fourier representation only for the efficient inverse of the preconditioner $\Kref$ as discussed in Section~\ref{sec:preconditioning}.
    
    \item \emph{Preconditioner and reference material.} 
Our preconditioner \eqref{eq:Kref} has the same form as the fundamental solution $\MB{G}^{0}$ contained in the discretized periodic Green's operator $\MB{\Gamma}_{(i)}^0$ of FFT-$Q_{1}$ Hex scheme (equation $(16)$ of~\cite{schneider_fft-based_2017}), and the fundamental solution $\widehat{\phi}$ in FANS (equation $(49)$ of~\cite{Leuschner2018}), therefore we expect similar conditioning of all three schemes. 
However, we provide detailed insight from a linear algebra viewpoint. Direct correspondence between the reference material $\rmaterial_{(i)}$, material ${\material_{(i)}}(\vek{x})$ and the resulting eigenvalues renders the optimization of $\rmaterial_{(i)}$ more accessible.
 The closer the reference material is to the real material of the sample, the better conditioning the discretized problem has. Therefore, in contrast to \cite{ZeGeVoPeGe2017,DEGEUS2017412}, we recommend reassembling the preconditioner $\KrefF$ when the material tangent significantly changes in Newton's method.


     %

    \item \emph{Discretisation grid.} Both FFT-$Q_{1}$ Hex and FANS are developed for bi/trilinear FE basis and quadrilateral/hexahedral elements. Their authors mentioned a possible extension for more complex elements which we present in this paper. In addition, the discretisation grid of our method does not have to follow the pixel/voxel structure. We allow for an arbitrary space-filling pattern of elements to be used, recall the patterns in Fig.~\ref{fig:tiling}. Further extension of our formulation to FE with higher-order polynomial basis functions is therefore straightforward. 
    
    \item \emph{Computational complexity.}
    Computational complexity of FFT-based methods is governed by $\mathcal{O}(n\log n$) complexity of the FFT. However, in our scheme, we compute two FFTs on $d\Nn$ displacement fields of size ${\NI}$, instead of $\mandeld$ strain fields of size ${\NQ}$ in FFT-$Q_{1}$. Because the number of strain components $\mandeld$ exceeds the number of displacement components $d\Nn$ per stencil and the number of quadrature points $\NQ$ exceeds the number of discretisation nodes $\NI$, our method has smaller computational overhead than the DB methods that evaluate the gradient in the Fourier space and perform FFT on the strain-sized fields. For instance, in the case of trilinear hexahedral FEs with $8$ quadrature points per element, the saving factor of our method is $24$.  
    
%
     %
\end{itemize} 

 \section{Conclusions}\label{sec:conclusion}
In this paper, we present a novel and \emph{optimal} approach for computational homogenization of nonlinear micromechanical and thermal problems in periodic media. The efficiency is achieved due to a clever interplay between the PCG solver and the geometry and physical properties of the problem \cite{MalekStakos}. Standard FE discretisation on a regular grid is coupled with the Newton's method to handle the nonlinear system iteratively.
The linearised system is solved by the PCG method, which is enhanced with a preconditioner based on a discretized inverse (Green's) operator for a problem with spatially uniform reference material data. The proposed matrix-free method exhibits excellent convergence properties as the number of linear solver iterations is bounded independently of the number of discretisation nodes and shows mild phase-contrast sensitivity. Our main findings are summarized as follows:
\begin{itemize}
    \item The condition number associated with the preconditioned linear system decreases as the reference material data approaches the material data. Two-sided bounds for all eigenvalues of the preconditioned linear system are easily accessible and thus provide valuable insight into the choice of the reference material.
         
    \item The computational complexity is governed by the FFT algorithm applied to the displacement field. 
    The preconditioning operator is cheaply inverted and applied in Fourier space, while the gradient is evaluated through the convolution with a short kernel in the real space.
    
     \item  
     The FE bases produce oscillation-free stress and strain solution fields with marginal discretisation artifacts at the phase interfaces. Additional variability of discretisation patterns allows the reduction of mesh anisotropy and a more accurate representation of the geometry and the solution.
\end{itemize}

 In addition, the Galerkin nature of the FE method connected with the minimization of the related energy functional allows us to use a well-built theory on the FE method for error estimation, convergence analysis, and other useful tools. 
In the future, the extension of the equivalence of DB and SB schemes to a general reference material  and the fusion of low-rank tensor approximation technique of \citet{VONDREJC2020112890} with our FE scheme are of primal interest.


\section*{Acknowledgments}
ML, IP, and JZ acknowledge the support by the Center of Advanced Applied Sciences, the European Regional Development Fund (project No. CZ 02.1.01/0.0/0.0/16 019/0000778). ML additionally acknowledges support by the Czech Science Foundation (project No. 20-14736S) and the Student Grant Competition of the Czech Technical University in Prague (project No. SGS22/004/OHK1/1T/11). RJL and LP acknowledge support from the European Research Council (StG-757343).

\appendix
\section{Appendix}\label{sec:appendix}
\paragraph{Thermal conduction}\label{sec:appendix_termal}
The proposed preconditioned FE method can be used also for potential problems such as thermal conduction or electrostatics. From a mathematical viewpoint, these problems are described by a scalar elliptic partial differential equation. 

For the scalar thermal conduction problem, we split the overall temperature gradient $\grad \temperature: \puc \rightarrow \D{R}^{d}$ into an average temperature gradient $\macrostrain =\frac{1}{| \puc |}\int_{\puc} \grad \temperature(\vek{x}) \intd{\vek{x}}  \in \D{R}^{d}$ and a periodically fluctuating field $\grad \pertemp : \puc \rightarrow \D{R}^{d}$
\begin{align*}
    \grad \temperature
    =
    \macrostrain+\grad \pertemp \quad \text{for all }  \vek{x} \in \puc.
\end{align*}
Here, $\grad \pertemp$ denotes the temperature gradient, and the fluctuating temperature field  $\pertemp$ belongs to the space of admissible functions
$\testspace=\left\lbrace  \tilde{v} :\puc\rightarrow\D{R},\, \tilde{v} \text{ is } \puc \text{-periodic }\right\rbrace$.
The governing equation for $\grad \pertemp$ follows from the thermal equilibrium condition
 \begin{align*}
     -\grad \cdot \flux(\vek{x},\macrostrain+\grad \pertemp(\vek{x}))=0 \quad \text{for all }  \vek{x} \in \puc,
 \end{align*}
in which $\flux:\puc \times \D{R}^{d}\times \D{R}^{q}\rightarrow \D{R}^{d}$ is the flux field. As usual, the equilibrium equation is converted to the weak form
\begin{align}\label{eq:weak_form_scalar}
    \int_{\puc}
   \grad \tilde{v}(\vek{x})^{\Transpose}
    \flux(\vek{x},\macrostrain+\grad \pertemp(\vek{x}))
    \intd{\vek{x}}
    = 0  \quad \text{for all } \tilde{v} \in \testspace
\end{align}
that serves as the starting point for the FE method.
Following the discretisation scheme described in Section~\ref{sec:FEM}, the linearisation in Section~\ref{sec:linearisation} and preconditioning in Section~\ref{sec:preconditioning} leads to a well-conditioned linear system
    \begin{align}
        \underbrace{\MB{F}^{H}(\MB{F}\Dmat^{\Transpose} \MB{\qweights} \conductDrmaterial_{(i)} \Dmat\MB{F}^{H} )^{-1}\MB{F}}_{\Krefinv}
         \underbrace{
         \Dmat^{\Transpose}
          \MB{\qweights}
        \conductDmaterial_{(i)}
         \Dmat
         }_{\linsysMat_{(i)}}
         \delta\Dpertemp_{(i+1)}
         =&
         \underbrace{
         \MB{F}^{H}(\MB{F}\Dmat^{\Transpose}  \MB{\qweights} \conductDrmaterial_{(i)} \Dmat\MB{F}^{H} )^{-1}\MB{F}
         }_{\Krefinv}
         \underbrace{
         \Dmat^{\Transpose} \MB{\qweights} \flux(\MB{e}+\Dmat\Dpertemp_{(i)})
         }_{-\RHS_{(i)}},
\end{align}
for a finite Newton's method increment $\delta\Dpertemp_{(i+1)}$. Material data matrix $\MB{\conductDmaterial}_{(i)}\in \D{R}^{d\NQ \times d\NQ}$ stores values of conductivity tangent matrix  $\MB{\conductmaterial}_{(i)}(\vek{x})=\dfrac{\partial \flux}{\partial\grad \pertemp}(\vek{x},\macrostrain+\grad \pertemp (\vek{x}) )\in \D{R}^{d\times d}$ in $(i)$-th Newton's method step, and $\conductDrmaterial_{(i)}\in \D{R}^{d\NQ \times d\NQ}$ comes from spatially uniform material data~$\conductrmaterial_{(i)}\in \D{R}^{d\times d}$. Another small difference lies in the form of the gradient matrix~$\Dmat$, 
    \begin{align}\label{eq:Dt_mat}
    \MB{\grad \Dpertemp}
    =
    \Dmat\MB{\Dpertemp}=
    {\begin{bmatrix}
   \Dmat_{1}  \\ 
   \Dmat_{2}  
   \end{bmatrix}}
    {\begin{bmatrix}
   \MB{\Dpertemp}
   \end{bmatrix}}.
    \end{align}
Here, the entries are the same as in the elasticity problem, recall equation~\eqref{eq:grad_componets}.
\bibliographystyle{elsarticle-num-names} 
\bibliography{cas-refs}

\begin{thebibliography}{41}
\expandafter\ifx\csname natexlab\endcsname\relax\def\natexlab#1{#1}\fi
\providecommand{\url}[1]{\texttt{#1}}
\providecommand{\href}[2]{#2}
\providecommand{\path}[1]{#1}
\providecommand{\DOIprefix}{doi:}
\providecommand{\ArXivprefix}{arXiv:}
\providecommand{\URLprefix}{URL: }
\providecommand{\Pubmedprefix}{pmid:}
\providecommand{\doi}[1]{\href{http://dx.doi.org/#1}{\path{#1}}}
\providecommand{\Pubmed}[1]{\href{pmid:#1}{\path{#1}}}
\providecommand{\bibinfo}[2]{#2}
\ifx\xfnm\relax \def\xfnm[#1]{\unskip,\space#1}\fi
\bibitem[{LLorca et~al.(2011)LLorca, González, Molina-Aldareguía, Segurado,
  Seltzer, Sket, Rodríguez, Sádaba, Muñoz, and Canal}]{LLorca_2011}
\bibinfo{author}{J.~LLorca}, \bibinfo{author}{C.~González},
  \bibinfo{author}{J.~M. Molina-Aldareguía}, \bibinfo{author}{J.~Segurado},
  \bibinfo{author}{R.~Seltzer}, \bibinfo{author}{F.~Sket},
  \bibinfo{author}{M.~Rodríguez}, \bibinfo{author}{S.~Sádaba},
  \bibinfo{author}{R.~Muñoz}, \bibinfo{author}{L.~P. Canal},
\newblock \bibinfo{title}{Multiscale modeling of composite materials: a roadmap
  towards virtual testing},
\newblock \bibinfo{journal}{Advanced Materials} \bibinfo{volume}{23}
  (\bibinfo{year}{2011}) \bibinfo{pages}{5130--5147}.
  \DOIprefix\doi{10.1002/adma.201101683}.
\bibitem[{Matouš et~al.(2017)Matouš, Geers, Kouznetsova, and
  Gillman}]{MATOUS2017192}
\bibinfo{author}{K.~Matouš}, \bibinfo{author}{M.~G.~D. Geers},
  \bibinfo{author}{V.~G. Kouznetsova}, \bibinfo{author}{A.~Gillman},
\newblock \bibinfo{title}{{A review of predictive nonlinear theories for
  multiscale modeling of heterogeneous materials}},
\newblock \bibinfo{journal}{Journal of Computational Physics}
  \bibinfo{volume}{330} (\bibinfo{year}{2017}) \bibinfo{pages}{192--220}.
  \DOIprefix\doi{10.1016/j.jcp.2016.10.070}.
\bibitem[{Fish et~al.(2021)Fish, Wagner, and Keten}]{Fish2021}
\bibinfo{author}{J.~Fish}, \bibinfo{author}{G.~J. Wagner},
  \bibinfo{author}{S.~Keten},
\newblock \bibinfo{title}{{Mesoscopic and multiscale modelling in materials}},
\newblock \bibinfo{journal}{Nature Materials} \bibinfo{volume}{20}
  (\bibinfo{year}{2021}) \bibinfo{pages}{774--786}.
  \DOIprefix\doi{10.1038/s41563-020-00913-0}.
\bibitem[{Terada et~al.(1997)Terada, Miura, and Kikuchi}]{Terada1997}
\bibinfo{author}{K.~Terada}, \bibinfo{author}{T.~Miura},
  \bibinfo{author}{N.~Kikuchi},
\newblock \bibinfo{title}{{Digital image-based modeling applied to the
  homogenization analysis of composite materials}},
\newblock \bibinfo{journal}{Computational Mechanics} \bibinfo{volume}{20}
  (\bibinfo{year}{1997}) \bibinfo{pages}{331--346}.
  \DOIprefix\doi{10.1007/s004660050255}.
\bibitem[{Maire and Withers(2014)}]{Maire2014}
\bibinfo{author}{E.~Maire}, \bibinfo{author}{P.~J. Withers},
\newblock \bibinfo{title}{Quantitative x-ray tomography},
\newblock \bibinfo{journal}{International Materials Reviews}
  \bibinfo{volume}{59} (\bibinfo{year}{2014}) \bibinfo{pages}{1--43}.
  \DOIprefix\doi{10.1179/1743280413Y.0000000023}.
\bibitem[{Sonon et~al.(2021)Sonon, {Ehab Moustafa Kamel}, and
  Massart}]{SONON20211}
\bibinfo{author}{B.~Sonon}, \bibinfo{author}{K.~{Ehab Moustafa Kamel}},
  \bibinfo{author}{T.~J. Massart},
\newblock \bibinfo{title}{{Advanced geometry representations and tools for
  microstructural and multiscale modeling}},
\newblock volume~\bibinfo{volume}{54} of \textit{\bibinfo{series}{Advances in
  Applied Mechanics}}, \bibinfo{publisher}{Elsevier}, \bibinfo{year}{2021}, pp.
  \bibinfo{pages}{1--111}. \DOIprefix\doi{10.1016/bs.aams.2020.12.001}.
\bibitem[{Johnson(1995)}]{johnson_1995}
\bibinfo{author}{C.~Johnson}, \bibinfo{title}{{Numerical Solution of Partial
  Differential Equations by the Finite Element Method}},
  \bibinfo{publisher}{Cambridge University Press}, \bibinfo{year}{1995}.
\bibitem[{Hestenes and Stiefel(1952)}]{Hestenes1952MethodsOC}
\bibinfo{author}{M.~R. Hestenes}, \bibinfo{author}{E.~Stiefel},
\newblock \bibinfo{title}{Methods of conjugate gradients for solving linear
  systems},
\newblock \bibinfo{journal}{Journal of research of the National Bureau of
  Standards} \bibinfo{volume}{49} (\bibinfo{year}{1952})
  \bibinfo{pages}{409--435}. \DOIprefix\doi{10.6028/JRES.049.044}.
\bibitem[{Moulinec and Suquet(1994)}]{moulinec_fast_1994}
\bibinfo{author}{H.~Moulinec}, \bibinfo{author}{P.~Suquet},
\newblock \bibinfo{title}{{A fast numerical method for computing the linear and
  nonlinear mechanical properties of composites}},
\newblock \bibinfo{journal}{{Comptes Rendus de l'Acad{\'e}mie des sciences.
  S{\'e}rie II. M{\'e}canique, physique, chimie, astronomie}}
  \bibinfo{volume}{318} (\bibinfo{year}{1994}) \bibinfo{pages}{1417--1423}.
  \URLprefix \url{https://hal.archives-ouvertes.fr/hal-03019226}.
\bibitem[{Moulinec and Suquet(1998)}]{moulinec_numerical_1998}
\bibinfo{author}{H.~Moulinec}, \bibinfo{author}{P.~Suquet},
\newblock \bibinfo{title}{{A numerical method for computing the overall
  response of nonlinear composites with complex microstructure}},
\newblock \bibinfo{journal}{Computer Methods in Applied Mechanics and
  Engineering} \bibinfo{volume}{157} (\bibinfo{year}{1998})
  \bibinfo{pages}{69--94}. \DOIprefix\doi{10.1016/S0045-7825(97)00218-1}.
\bibitem[{Golub and Van~Loan(2013)}]{golub2013matrix}
\bibinfo{author}{G.~Golub}, \bibinfo{author}{C.~Van~Loan},
  \bibinfo{title}{{Matrix Computations}}, Johns Hopkins Studies in the
  Mathematical Sciences, \bibinfo{publisher}{Johns Hopkins University Press},
  \bibinfo{year}{2013}.
\bibitem[{Schneider(2021)}]{Schneider2021}
\bibinfo{author}{M.~Schneider},
\newblock \bibinfo{title}{{A review of nonlinear FFT-based computational
  homogenization methods}},
\newblock \bibinfo{journal}{Acta Mechanica} \bibinfo{volume}{232}
  (\bibinfo{year}{2021}) \bibinfo{pages}{2051--2100}.
  \DOIprefix\doi{10.1007/s00707-021-02962-1}.
\bibitem[{Lucarini et~al.(2021)Lucarini, Upadhyay, and
  Segurado}]{Lucarini_2021}
\bibinfo{author}{S.~Lucarini}, \bibinfo{author}{M.~V. Upadhyay},
  \bibinfo{author}{J.~Segurado},
\newblock \bibinfo{title}{{FFT} based approaches in micromechanics:
  fundamentals, methods and applications},
\newblock \bibinfo{journal}{Modelling and Simulation in Materials Science and
  Engineering} \bibinfo{volume}{30} (\bibinfo{year}{2021})
  \bibinfo{pages}{023002}. \DOIprefix\doi{10.1088/1361-651x/ac34e1}.
\bibitem[{Brisard and Dormieux(2010)}]{Brisard2010FFT}
\bibinfo{author}{S.~Brisard}, \bibinfo{author}{L.~Dormieux},
\newblock \bibinfo{title}{{FFT-based methods for the mechanics of composites: A
  general variational framework}},
\newblock \bibinfo{journal}{Computational Materials Science}
  \bibinfo{volume}{49} (\bibinfo{year}{2010}) \bibinfo{pages}{663--671}.
  \DOIprefix\doi{10.1016/j.commatsci.2010.06.009}.
\bibitem[{Zeman et~al.(2010)Zeman, Vondřejc, Nov{\'{a}}k, and
  Marek}]{ZeVoNoMa2010AFFTH}
\bibinfo{author}{J.~Zeman}, \bibinfo{author}{J.~Vondřejc},
  \bibinfo{author}{J.~Nov{\'{a}}k}, \bibinfo{author}{I.~Marek},
\newblock \bibinfo{title}{{Accelerating a FFT-based solver for numerical
  homogenization of periodic media by conjugate gradients}},
\newblock \bibinfo{journal}{Journal of Computational Physics}
  \bibinfo{volume}{229} (\bibinfo{year}{2010}) \bibinfo{pages}{8065--8071}.
  \DOIprefix\doi{10.1016/j.jcp.2010.07.010}.
\bibitem[{Brisard and Dormieux(2012)}]{BRISARD2012197}
\bibinfo{author}{S.~Brisard}, \bibinfo{author}{L.~Dormieux},
\newblock \bibinfo{title}{{Combining Galerkin approximation techniques with the
  principle of Hashin and Shtrikman to derive a new FFT-based numerical method
  for the homogenization of composites}},
\newblock \bibinfo{journal}{Computer Methods in Applied Mechanics and
  Engineering} \bibinfo{volume}{217-220} (\bibinfo{year}{2012})
  \bibinfo{pages}{197--212}. \DOIprefix\doi{10.1016/j.cma.2012.01.003}.
\bibitem[{Vondřejc et~al.(2014)Vondřejc, Zeman, and Marek}]{VoZeMa2014FFTH}
\bibinfo{author}{J.~Vondřejc}, \bibinfo{author}{J.~Zeman},
  \bibinfo{author}{I.~Marek},
\newblock \bibinfo{title}{{An FFT-based Galerkin method for homogenization of
  periodic media}},
\newblock \bibinfo{journal}{Computers \& Mathematics with Applications}
  \bibinfo{volume}{68} (\bibinfo{year}{2014}) \bibinfo{pages}{156--173}.
  \DOIprefix\doi{10.1016/j.camwa.2014.05.014}.
\bibitem[{G{\'{e}}l{\'{e}}bart and Mondon-Cancel(2013)}]{Gelebart2013}
\bibinfo{author}{L.~G{\'{e}}l{\'{e}}bart}, \bibinfo{author}{R.~Mondon-Cancel},
\newblock \bibinfo{title}{{Non-linear extension of FFT-based methods
  accelerated by conjugate gradients to evaluate the mechanical behavior of
  composite materials}},
\newblock \bibinfo{journal}{Computational Materials Science}
  \bibinfo{volume}{77} (\bibinfo{year}{2013}) \bibinfo{pages}{430--439}.
  \DOIprefix\doi{10.1016/j.commatsci.2013.04.046}.
\bibitem[{Kabel et~al.(2014)Kabel, B{\"{o}}hlke, and
  Schneider}]{Kabel2014LargeDef}
\bibinfo{author}{M.~Kabel}, \bibinfo{author}{T.~B{\"{o}}hlke},
  \bibinfo{author}{M.~Schneider},
\newblock \bibinfo{title}{{Efficient fixed point and Newton–Krylov solvers
  for FFT-based homogenization of elasticity at large deformations}},
\newblock \bibinfo{journal}{Computational Mechanics} \bibinfo{volume}{54}
  (\bibinfo{year}{2014}) \bibinfo{pages}{1497--1514}.
  \DOIprefix\doi{10.1007/s00466-014-1071-8}.
\bibitem[{Zeman et~al.(2017)Zeman, de~Geus, Vondřejc, Peerlings, and
  Geers}]{ZeGeVoPeGe2017}
\bibinfo{author}{J.~Zeman}, \bibinfo{author}{T.~W.~J. de~Geus},
  \bibinfo{author}{J.~Vondřejc}, \bibinfo{author}{R.~H.~J. Peerlings},
  \bibinfo{author}{M.~G.~D. Geers},
\newblock \bibinfo{title}{{A finite element perspective on nonlinear FFT-based
  micromechanical simulations}},
\newblock \bibinfo{journal}{International Journal for Numerical Methods in
  Engineering} \bibinfo{volume}{111} (\bibinfo{year}{2017})
  \bibinfo{pages}{903--926}. \DOIprefix\doi{10.1002/nme.5481}.
\bibitem[{{de Geus} et~al.(2017){de Geus}, Vondřejc, Zeman, Peerlings, and
  Geers}]{DEGEUS2017412}
\bibinfo{author}{T.~{de Geus}}, \bibinfo{author}{J.~Vondřejc},
  \bibinfo{author}{J.~Zeman}, \bibinfo{author}{R.~Peerlings},
  \bibinfo{author}{M.~Geers},
\newblock \bibinfo{title}{{Finite strain FFT-based non-linear solvers made
  simple}},
\newblock \bibinfo{journal}{Computer Methods in Applied Mechanics and
  Engineering} \bibinfo{volume}{318} (\bibinfo{year}{2017})
  \bibinfo{pages}{412--430}. \DOIprefix\doi{10.1016/j.cma.2016.12.032}.
\bibitem[{Kaßbohm et~al.(2006)Kaßbohm, Müller, and Feßler}]{Kasbohm2006}
\bibinfo{author}{S.~Kaßbohm}, \bibinfo{author}{W.~H. Müller},
  \bibinfo{author}{R.~Feßler},
\newblock \bibinfo{title}{Improved approximations of fourier coefficients for
  computing periodic structures with arbitrary stiffness distribution},
\newblock \bibinfo{journal}{Computational Materials Science}
  \bibinfo{volume}{37} (\bibinfo{year}{2006}) \bibinfo{pages}{90--93}.
  \DOIprefix\doi{10.1016/j.commatsci.2005.12.010}.
\bibitem[{Shanthraj et~al.(2015)Shanthraj, Eisenlohr, Diehl, and
  Roters}]{SHANTHRAJ201531}
\bibinfo{author}{P.~Shanthraj}, \bibinfo{author}{P.~Eisenlohr},
  \bibinfo{author}{M.~Diehl}, \bibinfo{author}{F.~Roters},
\newblock \bibinfo{title}{Numerically robust spectral methods for crystal
  plasticity simulations of heterogeneous materials},
\newblock \bibinfo{journal}{International Journal of Plasticity}
  \bibinfo{volume}{66} (\bibinfo{year}{2015}) \bibinfo{pages}{31--45}.
  \DOIprefix\doi{10.1016/j.ijplas.2014.02.006}.
\bibitem[{Willot et~al.(2014)Willot, Abdallah, and Pellegrini}]{Willot2014}
\bibinfo{author}{F.~Willot}, \bibinfo{author}{B.~Abdallah},
  \bibinfo{author}{Y.-P. Pellegrini},
\newblock \bibinfo{title}{{Fourier-based schemes with modified Green operator
  for computing the electrical response of heterogeneous media with accurate
  local fields}},
\newblock \bibinfo{journal}{International Journal for Numerical Methods in
  Engineering} \bibinfo{volume}{98} (\bibinfo{year}{2014})
  \bibinfo{pages}{518--533}. \DOIprefix\doi{10.1002/nme.4641}.
\bibitem[{Schneider et~al.(2016)Schneider, Ospald, and Kabel}]{Schneider2016}
\bibinfo{author}{M.~Schneider}, \bibinfo{author}{F.~Ospald},
  \bibinfo{author}{M.~Kabel},
\newblock \bibinfo{title}{{Computational homogenization of elasticity on a
  staggered grid}},
\newblock \bibinfo{journal}{International Journal for Numerical Methods in
  Engineering} \bibinfo{volume}{105} (\bibinfo{year}{2016})
  \bibinfo{pages}{693--720}. \DOIprefix\doi{10.1002/nme.5008}.
\bibitem[{Schneider et~al.(2017)Schneider, Merkert, and
  Kabel}]{schneider_fft-based_2017}
\bibinfo{author}{M.~Schneider}, \bibinfo{author}{D.~Merkert},
  \bibinfo{author}{M.~Kabel},
\newblock \bibinfo{title}{{FFT-based homogenization for microstructures
  discretized by linear hexahedral elements}},
\newblock \bibinfo{journal}{International Journal for Numerical Methods in
  Engineering} \bibinfo{volume}{109} (\bibinfo{year}{2017})
  \bibinfo{pages}{1461--1489}. \DOIprefix\doi{10.1002/nme.5336}.
\bibitem[{Leuschner and Fritzen(2018)}]{Leuschner2018}
\bibinfo{author}{M.~Leuschner}, \bibinfo{author}{F.~Fritzen},
\newblock \bibinfo{title}{{Fourier-Accelerated Nodal Solvers (FANS) for
  homogenization problems}},
\newblock \bibinfo{journal}{Computational Mechanics} \bibinfo{volume}{62}
  (\bibinfo{year}{2018}) \bibinfo{pages}{359--392}.
  \DOIprefix\doi{10.1007/s00466-017-1501-5}.
\bibitem[{Leute et~al.(2022)Leute, Ladecký, Falsafi, Jödicke, Pultarová,
  Zeman, Junge, and Pastewka}]{LeuteR2021}
\bibinfo{author}{R.~J. Leute}, \bibinfo{author}{M.~Ladecký},
  \bibinfo{author}{A.~Falsafi}, \bibinfo{author}{I.~Jödicke},
  \bibinfo{author}{I.~Pultarová}, \bibinfo{author}{J.~Zeman},
  \bibinfo{author}{T.~Junge}, \bibinfo{author}{L.~Pastewka},
\newblock \bibinfo{title}{{Elimination of ringing artifacts by finite-element
  projection in FFT-based homogenization}},
\newblock \bibinfo{journal}{Journal of Computational Physics}
  \bibinfo{volume}{453} (\bibinfo{year}{2022}) \bibinfo{pages}{110931}.
  \DOIprefix\doi{10.1016/j.jcp.2021.110931}.
\bibitem[{Ma et~al.(2021)Ma, Shakoor, Vasiukov, Lomov, and Park}]{Ma2021}
\bibinfo{author}{X.~Ma}, \bibinfo{author}{M.~Shakoor},
  \bibinfo{author}{D.~Vasiukov}, \bibinfo{author}{S.~V. Lomov},
  \bibinfo{author}{C.~H. Park},
\newblock \bibinfo{title}{{Numerical artifacts of Fast Fourier Transform
  solvers for elastic problems of multi-phase materials: their causes and
  reduction methods}},
\newblock \bibinfo{journal}{Computational Mechanics}  (\bibinfo{year}{2021}).
  \DOIprefix\doi{10.1007/s00466-021-02013-5}.
\bibitem[{Axelsson and Karátson(2009)}]{Axelsson2009}
\bibinfo{author}{O.~Axelsson}, \bibinfo{author}{J.~Karátson},
\newblock \bibinfo{title}{Equivalent operator preconditioning for elliptic
  problems},
\newblock \bibinfo{journal}{Numerical Algorithms} \bibinfo{volume}{50}
  (\bibinfo{year}{2009}) \bibinfo{pages}{297--380}.
  \DOIprefix\doi{10.1007/s11075-008-9233-4}.
\bibitem[{Pultarová and Ladecký(2021)}]{PultarovaNLAA}
\bibinfo{author}{I.~Pultarová}, \bibinfo{author}{M.~Ladecký},
\newblock \bibinfo{title}{Two-sided guaranteed bounds to individual
  eigenvalues of preconditioned finite element and finite difference
  problems},
\newblock \bibinfo{journal}{Numer Linear Algebra Appl.} \bibinfo{volume}{28}
  (\bibinfo{year}{2021}) \bibinfo{pages}{e2382}.
  \DOIprefix\doi{https://doi.org/10.1002/nla.2382}.
\bibitem[{Gergelits et~al.(2019)Gergelits, Mardal, Nielsen, and
  Strakoš}]{Gergelits_2019}
\bibinfo{author}{T.~Gergelits}, \bibinfo{author}{K.-A. Mardal},
  \bibinfo{author}{B.~F. Nielsen}, \bibinfo{author}{Z.~Strakoš},
\newblock \bibinfo{title}{{Laplacian Preconditioning of Elliptic PDEs:
  Localization of the Eigenvalues of the Discretized Operator}},
\newblock \bibinfo{journal}{SIAM Journal on Numerical Analysis}
  \bibinfo{volume}{57} (\bibinfo{year}{2019}) \bibinfo{pages}{1369--1394}.
  \DOIprefix\doi{10.1137/18M1212458}.
\bibitem[{Ladecký et~al.(2020)Ladecký, Pultarová, and Zeman}]{Ladecky2020}
\bibinfo{author}{M.~Ladecký}, \bibinfo{author}{I.~Pultarová},
  \bibinfo{author}{J.~Zeman},
\newblock \bibinfo{title}{{Guaranteed Two-Sided Bounds on All Eigenvalues of
  Preconditioned Diffusion and Elasticity Problems Solved By the Finite Element
  Method}},
\newblock \bibinfo{journal}{Applications of Mathematics} \bibinfo{volume}{66}
  (\bibinfo{year}{2020}) \bibinfo{pages}{21--42}.
  \DOIprefix\doi{10.21136/AM.2020.0217-19}.
\bibitem[{Cooley and Tukey(1965)}]{CooleyTukey}
\bibinfo{author}{J.~W. Cooley}, \bibinfo{author}{J.~W. Tukey},
\newblock \bibinfo{title}{{An algorithm for the machine calculation of complex
  Fourier series}},
\newblock \bibinfo{journal}{Mathematics of Computation} \bibinfo{volume}{19}
  (\bibinfo{year}{1965}) \bibinfo{pages}{297--301}.
  \DOIprefix\doi{10.1090/S0025-5718-1965-0178586-1}.
\bibitem[{Saad(2003)}]{Saad2003}
\bibinfo{author}{Y.~Saad}, \bibinfo{title}{Iterative Methods for Sparse Linear
  Systems}, \bibinfo{edition}{second} ed., \bibinfo{publisher}{Society for
  Industrial and Applied Mathematics}, \bibinfo{year}{2003}.
  \DOIprefix\doi{10.1137/1.9780898718003}.
\bibitem[{Eijkhout and Vassilevski(1991)}]{Eijkhout1991TheRO}
\bibinfo{author}{V.~Eijkhout}, \bibinfo{author}{P.~Vassilevski},
\newblock \bibinfo{title}{{The Role of the Strengthened
  Cauchy-Buniakowskii-Schwarz Inequality in Multilevel Methods}},
\newblock \bibinfo{journal}{SIAM Rev.} \bibinfo{volume}{33}
  (\bibinfo{year}{1991}) \bibinfo{pages}{405--419}.
  \DOIprefix\doi{10.1137/1033098}.
\bibitem[{Axelsson(1996)}]{axelsson_1996}
\bibinfo{author}{O.~Axelsson}, \bibinfo{title}{Iterative Solution Methods},
  \bibinfo{publisher}{Cambridge University Press}, \bibinfo{year}{1996}.
  \DOIprefix\doi{10.1017/CBO9780511624100}.
\bibitem[{Nielsen et~al.(2009)Nielsen, Tveito, and
  Hackbusch}]{nielsen_preconditioning_2009}
\bibinfo{author}{B.~F. Nielsen}, \bibinfo{author}{A.~Tveito},
  \bibinfo{author}{W.~Hackbusch},
\newblock \bibinfo{title}{Preconditioning by inverting the {Laplacian}: an
  analysis of the eigenvalues},
\newblock \bibinfo{journal}{IMA Journal of Numerical Analysis}
  \bibinfo{volume}{29} (\bibinfo{year}{2009}) \bibinfo{pages}{24--42}.
  \DOIprefix\doi{10.1093/imanum/drm018}.
\bibitem[{Ladecký et~al.(tion)Ladecký, Falsafi, Leute, Pultarová, Zeman,
  Junge, and Pastewka}]{Ladecky2022}
\bibinfo{author}{M.~Ladecký}, \bibinfo{author}{A.~Falsafi},
  \bibinfo{author}{R.~J. Leute}, \bibinfo{author}{I.~Pultarová},
  \bibinfo{author}{J.~Zeman}, \bibinfo{author}{T.~Junge},
  \bibinfo{author}{L.~Pastewka},
\newblock \bibinfo{title}{{On the equivalence of the displacement- and
  strain-based computational homogenisation schemes}}  (\bibinfo{year}{in
  preparation}).
\bibitem[{Málek and Strakoš(2015)}]{MalekStakos}
\bibinfo{author}{J.~Málek}, \bibinfo{author}{Z.~Strakoš},
  \bibinfo{title}{{Preconditioning and the Conjugate Gradient Method in the
  Context of Solving PDEs}}, SIAM Spotlight Series, \bibinfo{publisher}{Society
  for Industrial and Applied Mathematics}, \bibinfo{year}{2015}.
\bibitem[{Vondřejc et~al.(2020)Vondřejc, Liu, Ladecký, and
  Matthies}]{VONDREJC2020112890}
\bibinfo{author}{J.~Vondřejc}, \bibinfo{author}{D.~Liu},
  \bibinfo{author}{M.~Ladecký}, \bibinfo{author}{H.~G. Matthies},
\newblock \bibinfo{title}{{FFT-based homogenisation accelerated by low-rank
  tensor approximations}},
\newblock \bibinfo{journal}{Computer Methods in Applied Mechanics and
  Engineering} \bibinfo{volume}{364} (\bibinfo{year}{2020})
  \bibinfo{pages}{112890}. \DOIprefix\doi{10.1016/j.cma.2020.112890}.

\end{thebibliography}





\end{document}